\documentclass[11pt]{amsart}
\usepackage{epsfig}
\usepackage{graphicx, subfigure}
\usepackage{amsmath,amsthm,amssymb,amsfonts,mathtools,mathrsfs,enumerate,mathabx}

\newtheorem{thm}{Theorem}[section]

\newtheorem{lem}[thm]{Lemma}
\newtheorem{prop}[thm]{Proposition}
\theoremstyle{definition}
\newtheorem{defn}[thm]{Definition}
\theoremstyle{remark}
\newtheorem{rem}[thm]{Remark}

\numberwithin{equation}{section}

\newcommand{\eps}{\varepsilon}

\begin{document}


\title[Global diffusion on a tight three-sphere]{Global diffusion on a tight three-sphere}


\author[Marian Gidea]{Marian Gidea$^\ddag$}
\address{Department of Mathematical Science, Yeshiva University, NY 10018, USA}%
\email{Marian.Gidea@yu.edu}%
\thanks{$^\ddag$ Research of M.G. was partially supported by NSF grant  DMS-0635607 and DMS-1201357.}

\begin{abstract}
We consider  an integrable Hamiltonian system weakly coupled with a pendulum-type system. For each energy level within some range, the uncoupled system   is assumed to possess a normally hyperbolic invariant manifold diffeomorphic to a three-sphere, which  bounds a strictly convex domain, and whose  stable and unstable invariant manifolds  coincide. The Hamiltonian flow on the three-sphere is equivalent to the  Reeb flow for the induced contact form. The strict convexity condition implies that the contact structure on the three-sphere is tight. When a  small, generic coupling is added to the system,   the normally hyperbolic invariant manifold is preserved as a three-sphere, and the stable and unstable manifolds split, yielding transverse intersections.
We show that there exist trajectories that follow any prescribed  collection  of invariant tori and Aubry-Mather sets  within some global section of the flow restricted to the three-sphere.
In this sense, we say that the perturbed system exhibits global diffusion on  the tight three-sphere.
\end{abstract}


\maketitle

\section{Introduction}

The diffusion problem in Hamiltonian dynamics asserts that,  generically, nearly integrable Hamiltonian systems exhibit  trajectories along which the action variable changes by some positive distance that is independent of the size of the perturbation (see \cite{Arnold64}). A related phenomenon is the existence of trajectories that exhibit symbolic dynamics, i.e., they visit prescribed open sets in the action variable domain, or prescribed sequences of invariant objects (KAM tori, Aubry-Mather sets, etc.) in the phase space. In this paper we are concerned with the phenomenon of diffusion for a particular  subclass of nearly integrable Hamiltonian systems,  referred to as {\em a priori unstable systems} -- these can be locally described in terms of action-angle and hyperbolic variables  (see \cite{ChierchiaG94}).
The diffusion phenomenon in a priori unstable systems has been extensively studied,   for instance in   \cite{ChierchiaG94,Treschev02b,BessiCV01,Berti2003,Kaloshin2003,Treschev04,
ChengY2004,LochakM05,DelshamsLS2006,GideaL06,Bernard08,ChengY2009,delshams2009geography,GideaR12}.

A typical approach to the diffusion problem is to start with an action-angle domain of the integrable Hamiltonian, which is foliated by  Liouville tori, apply the KAM theorem to the perturbed system restricted to that region, and use  KAM tori in combination with other geometric objects created by the perturbation to construct diffusing trajectories. Such trajectories lie within the action-angle domain which was originally populated by Liouville tori. In this sense, this type of diffusion is a {\em local phenomenon}.

The {\em global geometry} of an integrable system can, however, be quite rich, with the phase space divided out by the singular  leaves of the Liouville foliation into multiple action-angle domains, which are populated by different families of Liouville tori. In this paper we formulate a question on whether, under some appropriate conditions, small perturbations of such systems exhibit diffusing orbits that  travel arbitrarily far and across such multiple domains, and follow  some prescribed collections of invariant objects within each of these domains. We shall refer to this phenomenon as {\em global diffusion}. (We note that a different type of phenomenon, also referred to as `global diffusion', is studied in
\cite{Guzzo_firstnumerical}).

We investigate the question on global diffusion in a model  of an a priori unstable Hamiltonian system  of three degrees of freedom. The unperturbed  Hamiltonian system is a product of a two-degree  of freedom integrable Hamiltonian and a one-degree of freedom pendulum. We  assume that there exists a normally hyperbolic invariant manifold that is diffeomorphic to the product of a  three-sphere and an energy interval -- each three-sphere is normally hyperbolic within the corresponding energy level.  Furthermore, the two-degree of freedom integrable Hamiltonian is assumed to satisfy a strict convexity condition. Under the integrability condition that we assume,  each three-sphere is divided out by  singular leaves
into finitely many disjoint, open domains, each  described by one action and two angle coordinates. The stable and unstable invariant manifolds of each of the three-spheres coincide. We apply a small perturbation. The normally hyperbolic invariant manifold survives if the perturbation is small enough, and  its intersections with energy levels of the perturbed Hamiltonian yields a family of three-spheres, each being normally hyperbolic within the corresponding energy level.
Assuming  some non-degeneracy conditions that are generic,  the stable and unstable invariant manifolds intersect transversally along homoclinic manifolds.  To each homoclinic manifold, we associate a scattering map that describes the asymptotic behavior of the homoclinic trajectories.  Each scattering map is defined on some open subset of the normally hyperbolic invariant manifold. Typically, there are many geometrically distinct  homoclinic manifolds and corresponding scattering maps. We assume that there exists a family of homoclinic manifolds with the property that the domains of the corresponding scattering maps sweep across all possible action coordinate level sets on each action-angle domain  on the three-sphere. Under these assumptions, we show that there exist trajectories that follow  prescribed collections of invariant objects of a certain type within the three-sphere. In this sense, we say  that the perturbed system exhibits global diffusion relative to the three-sphere.

There are two key ingredients that are  essential for  us to formulate our global diffusion type of result.

The first ingredient is the assumption  that the domains of the scattering maps cover all level sets of the action coordinate on each action-angle domain on the three-sphere. This condition is satisfied  in very general situations.

The second ingredient is the strict convexity condition on the above mentioned three-spheres. A deep result from \cite{HZW98} implies that for each such a three-sphere, there exists  a disk-like global surface of section for the flow restricted to that sphere, and the Poincar\'e return map to the surface of section is equivalent to an area preserving map. We show that, when the two-degree of freedom  Hamiltonian is coupled with the pendulum via a small perturbation, as these  three-spheres survive as normally hyperbolic invariant manifolds in the corresponding energy levels, there still exists a global surface of section for the flow restricted to each three-sphere. We use  these global surfaces of section to describe the existence of trajectories that exhibit global diffusion. Specifically, within each action-angle domain determined within a global surface of section  we prescribe a collection of invariant objects that are either KAM invariant tori or Aubry-Mather sets. We show that  there exist   trajectories with visits each of these objects in the prescribed order. Again, the main novelty here, in comparison to other related works, is that {\em these invariant objects lie in different   action-angle domains}.

A motivation for the question on global diffusion is furnished  by the dynamics of the spatial circular restricted three-body problem (SCRTBP) near one of its equilibrium points. The SCRTBP  describes the spatial motion of an infinitesimal particle (e.g., a satellite) under the gravity of two heavy masses, referred to as primaries (e.g., Earth and Sun), moving on circular orbits about their common center of mass. It is assumed that the infinitesimal particle exerts no influence on the primaries. The motion of the infinitesimal particle relative to a co-rotating system of coordinates can be described by a Hamiltonian system of three degrees of freedom. One of the equilibrium points of the system is located between the primaries, and we refer to it as $L_1$. The linear stability of $L_1$ is of center$\times$center$\times$saddle  type. There exists a 4-dimensional center manifold through $L_1$ and a natural   action-angle  coordinate system  near $L_1$.  The intersections  of the  center manifold with energy levels above but close to that of $L_1$ are three-spheres. Each three-sphere is normally hyperbolic within its energy level, and   can be parametrized by one action and two angle coordinates, in such a way  that the action  corresponds to the out-of-plane amplitude of the motion. Provided that the mass ratio of the primaries is sufficiently small, and the energy level is sufficiently close to that of $L_1$, one can use  analytical results combined with   computer experiments  to show   that there exist trajectories that visit closely any prescribed sequence of action level sets -- at least within some bounds -- on the three-sphere.
Such an argument is carried  in \cite{Sama04}, and will be further developed in an upcoming work.
The same mechanism of global diffusion can be established via  numerical experiments for `realistic values' of mass ratios and energy levels; see  \cite{DelshamsGR10,DelshamsGR13}. Here by realistic values we mean those corresponding to
applications to astrodynamics.
A practical question motivating the global diffusion mechanism discussed in this paper is whether there exist trajectories that  visit any prescribed sequence of action level sets. This would
yield a zero cost procedure to design  satellite trajectories that change their out-of-plane amplitude, relative to the ecliptic, from nearly zero  to  nearly the maximum possible that can be achieved for a given energy level. Numerical evidence suggests that this is possible, at least for some convenient energy levels.

The model considered in the paper can be regarded as an {\em abstract model} for the SCRTBP,  for which a  purely analytical diffusion type of result can be established.
Similarities and differences between the Hamiltonian model considered in this paper and the  SCRTBP  are discussed in Remark \ref{rem:RTBP}.
Besides proposing the study of global diffusion,  in this paper we hope to further explore the integration of symplectic dynamics ideas into problems in Hamiltonian instability and celestial mechanics.  See, e.g.,   \cite{Albers11,Bramham11,Hry11,SalomaodePaulo13}.

Section \ref{section:first} provides background  and the main result. Section \ref{sec:application} discusses some examples. Section \ref{section:proofs} gives the proofs of the main result and of some auxiliary propositions.
The focus of the proof falls on describing the geometric structures that organize the dynamics. The actual argument for the existence of diffusing orbits follows closely the approach  in \cite{GideaR12}, therefore many details of the proof are omitted.
For easy reference, we are also including an Appendix.

\section{Set-up and main result}
\label{section:first}
We start by describing certain structures that are needed to formulate the main result of this paper. Then we describe a class of nearly integrable  Hamiltonian systems, and we formulate a global diffusion-type of result.

\subsection{Preliminaries}\subsubsection{Strictly convex energy hypersurfaces for two-degrees of freedom Hamiltonian systems}
Let $H:\mathbb{R}^4\to\mathbb{R}$  be a smooth Hamiltonian, where    $\mathbb{R}^4$ is endowed with the canonical symplectic form $\omega=\sum_{i=1}^{2}dq_i\wedge dp_i$. Let $\lambda=\frac{1}{2}\sum_{i=1}^{2}(q_idp_i-p_idq_i)$ be the Liouville  form.

Suppose that $c$ is a regular value for $H$ and let $\Lambda=\{(p_1,q_1,p_2,q_2)\,|\,H(p_1,q_1,p_2,q_2)=c\}$ be the corresponding energy hypersurface  of $H$. We say that $\Lambda$ is  star-shaped relative to a point $x_0$ if  \begin{itemize}
\item [(i)] for every $x\in\Lambda$, the ray $x_0+t(x-x_0)$, $t>0$,  intersects $\Lambda$ only once, at $x$, and
\item[(ii)] $(x-x_0)\not\in T_x\Lambda$.\end{itemize}
Then $\Lambda$ is diffeomorphic to  the $3$-dimensional unit sphere $S^3$ in $\mathbb{R}^4$, and there exists a positive smooth function $g:S^3\to\mathbb{R}^+$,  such that $\Lambda=\{z\sqrt{g(z)}\,|\,z\in S^3\}$ (see \cite{Eliashberg92}).

The restriction of $\lambda$ to $\Lambda$ is a contact form.  It determines a contact structure on $\Lambda$, which is the plane bundle  given by $\xi=\textrm{ker}(\lambda)$.

Contact forms can be overtwisted, if there exists an overtwisted disk, i.e., an embedded $2$-dimensional open disk $D$ with
$T\partial D\subseteq \xi$ and $T_pD\neq \xi_p$ for all $p\in\partial D$, and tight, otherwise.
Since $\Lambda$ is star-shaped, the contact form ${\lambda}_{\mid \Lambda}$ is equivalent to the contact form $g{\lambda}_{\mid S^3}$. This  contact form is tight.

The Reeb vector field $X_{\lambda}$ associated to ${\lambda}$ is defined by $i_{X_{\lambda}}(\lambda)=1$ and $i_{X_{\lambda}}(d\lambda )=0$, and the flow generated by the Reeb vector field is called the Reeb flow. A sufficient condition  for the Reeb flow to be conjugated to the Hamiltonian flow restricted to the energy hypersurface  is that the Liouville vector field $\eta$, defined by $L_{\eta}\omega=\omega$,  is transverse to $\Lambda$ (here $L$ denotes the Lie derivative).

The  Liouville vector field   is transverse to the star-shaped energy hypersurface, hence the Hamiltonian flow restricted to   $\Lambda$   is conjugate (via reparametrization of orbits) to the Reeb flow on $\Lambda$ defined
by the tight contact form ${g\lambda}_{\mid S^3}$.

An energy hypersurface $\Lambda$ is said to bound a strictly convex domain if there exists $C>0$ such that $(D^2H_{00} - C\cdot \textrm{id})$ is positive definite at all points on $\mathbb{R}^4$. If an energy hypersurface bounds a strictly convex domain, it implies that it is star-shaped with respect to some point.

An energy hypersurface  $\Lambda$ is said to be dynamically convex if every periodic solution $x$ of period $T$ of the Reeb flow on  $\Lambda$ has Conley-Zehnder index $\tilde\mu(x,T)\geq 3$. Dynamical convexity is a symplectic invariant. An energy hypersurface that bounds a strictly convex domain is always dynamically convex.
See Appendix \ref{section:symplectic} for further details.

\subsubsection{Integrable two-degrees of freedom Hamiltonian systems}\label{subsection:integrable}
Consider a Hamiltonian system $H:\mathbb{R}^4\to\mathbb{R}$.
The Hamiltonian system $H$ is said to be Liouville integrable if there exists a smooth function $K:\mathbb{R}^4\to\mathbb{R}$ such that
$H$ and $K$ are functionally  independent a.e. and in involution, i.e. $\{H,K\}=0$, where $\{\cdot,\cdot\}$ denotes the Poisson bracket.  The function $K$ is said to be a first integral of the Hamiltonian system. The mapping $\mathcal{H}:\mathbb{R}^4\to \mathbb{R}^2$, $\mathcal{H}=(H,K)$ is called the momentum mapping.

The momentum mapping determines a foliation of the phase space, called the Liouville foliation, given by the connected components of $\mathcal{L}=\{\mathcal{H}^{-1}(\bar a)\,|\, \bar a\in\mathbb{R}^2\}$. By the Liouville-Arnold Theorem \cite{Arnold63}, every compact and connected regular leaf of $\mathcal{L}$ is diffeomorphic to a two-torus $\mathbb{T}^2=\mathbb{R}^2/\mathbb{Z}^2$ and has a  neighborhood $U$, diffeomorphic to $D^2\times\mathbb{T}^2$, where $D^2$ is an open disk in $\mathbb{R}^2$, and there exists a canonical system of action-angle coordinates $(I_1,I_2,\phi_1,\phi_2)$ on $U$, $I=(I_1,I_2)\in D^2$, $\phi=(\phi_1,\phi_2)\in\mathbb{T}^2$, with $\sum_{i=1}^{2}dq_i\wedge dp_i=\sum_{i=1}^{2}d\phi_i\wedge dI_i$,  such that the Hamiltonian $H$  on $U$  depends only on $I$, and on each torus $T_{I}=\{I=\textrm{const.}\}$  the induced Hamiltonian flow is linear. The domain $U$ is called an action-angle domain. Each connected component $U$ of the complement of the set of the singular leaves of the Liouville foliation is a maximal action-angle domain.

Fix $\Lambda$   a regular energy hypersurface of $H$. We say that $K$ is a Bott function if the set of critical points of $K$ on $\Lambda$ is a disjoint union of smooth submanifolds, each
of which being non-degenerate in the following sense: $D^2K$ is
non-degenerate on the   transversals to the submanifold (at each point).\footnote{For the definition of a Bott function, it is not necessary to assume that  $H,K$ are functionally  independent a.e.    and in involution, as the  Bott integrability condition is a specific property of the restriction of
the Hamiltonian system to an energy level.}
The class of Bott systems was introduced  in \cite{Fomenko86} in the
study of the topology of integrable systems (see also \cite{BolsinovF04}).
It is known that the set of Bott systems on a given energy level is a set of first category in the set of all integrable
Hamiltonian systems with the weak $C^r$ metric, $r>2$ (see \cite{Kalashnikov2}).

The integral $K$ cannot have isolated critical points on $\Lambda$; any integral curve is non-degenerate and consists entirely of critical points, thus filling a submanifold of $\Lambda$.
It turns out that any connected critical submanifolds  of a Bott function on a compact, regular level set of a Hamiltonian of two degrees of freedom is diffeomorphic either to a circle, or to a torus, or to the Klein bottle. By an arbitrary small $C^r$-perturbation,  one can
turn an integrable Hamiltonian system with a Bott integral into an integrable
Hamiltonian system whose only critical submanifolds are circles (see \cite{Kalashnikov1,Kalashnikov2}).
If these circles are either elliptic or hyperbolic orbits for the Hamiltonian flow, the system is said to be coherent.

\subsubsection{Aubry-Mather theory}\label{subsection:aubrymather}
Assume that $H$ satisfies an iso-energetic non-degeneracy condition on a domain $D^2\times\mathbb{T}^2$ as above (see Section \ref{subsection:main} for explicit formulations of the iso-energetic non-degeneracy condition). Assume  that $\Sigma$ is a local surface of section for the Hamiltonian flow on $\Lambda$. Assume that $\mathcal{T}_1$ and $\mathcal{T}_2$ are two two-dimensional tori  in $\Lambda$ that are invariant under the flow. Assume that the region in $\Lambda$ between ${T}_1$ and $ {T}_2$ intersects $\Sigma$ in an annulus $A\simeq \mathbb{T}^1\times [0,1]$, bounded by two circles $ {\mathcal{T}}_1$ and $ {\mathcal{T}}_2$, and that the Poincar\'e first return map $f$ to $\Sigma$ is well defined on $A$.

Let $\mathbb{R}\times [0,1]$ be the universal cover, let $(\theta,J)\in\mathbb{R}\times [0,1]$ be a coordinate system on $A$, $\pi_J$ be the projection onto the first component, and $\pi_\theta$ be the projection onto the second component. We fix a lift of $f$ to $\mathbb{R}\times [0,1]$.  With an abuse of notation, we also denote the universal cover by $A$  and the lift by $f$.

In the setting of Section \ref{subsection:integrable}, the annulus $A$ is foliated by essential invariant circles. By an  essential invariant circle  we mean  a circle invariant under $f$  that cannot
be homotopically deformed into a point inside the annulus.
The iso-energetic non-degeneracy condition implies that $f$ satisfies a monotone twist-condition on $A$, that is,
 $\partial (\pi_{\theta}\circ  f)/\partial J\neq 0$ at all points in the annulus $A$.
Let us assume that   $f$ is a positive twist, meaning that $\partial (\pi_{\theta}\circ f)/\partial J>0$ at all points.

We remark here that the twist condition is a coordinate dependent condition.

Under some condition, we can obtain a perturbed  map on $A$, which we still denote by $f$, such that the  boundaries $ {\mathcal{T}}_1$ and $ {\mathcal{T}}_2$ of $A$ remain invariant under $f$, and $f$ is still an area preserving positive  twist map on $A$;   however  $A$ is no longer foliated by essential invariant circles. Such a situation  is described explicitly in Section \ref{subsection:main}.

The map $f$ restricted to the boundary components $ {\mathcal{T}}_1$, $ {\mathcal{T}}_2$ of the annulus  has well defined rotation numbers, which we denote $\rho_1,\rho_2$, respectively, with $\rho_1<\rho_2$.

An annulus $A$  as above is called a Birkhoff Zone of Instability (BZI) provided that there is no  essential invariant circle in the interior of the region.

A invariant subset $M\subseteq A$ is said to be monotone (cyclically ordered) if $\pi_\theta(z_1)<\pi_\theta(z_2)$ implies $\pi_\theta(f(z_1))<\pi_\theta(f(z_2))$  for all $z_1,z_2\in M$. For $z\in A$ the extended orbit of $z$ is the set $EO(z)=\{f^n(z)+(j,0)\,:\,n,j\in\mathbb{Z}\}$. The orbit of $z$ is said to be monotone  if the set $EO(z)$  is monotone. If the orbit of $z\in A$ is monotone, then its  rotation number $\rho(z)$ is well defined.
All points in the same monotone set have the same rotation number.

An Aubry-Mather set is a minimal, monotone,  $f$-invariant subset of given rotation number $\rho$.
Here by a minimal set we mean a closed invariant set that does not contain any proper closed invariant subsets. (Equivalently, the orbit of every point in the set is dense in the set.)  This should not be confused with action-minimizing sets.

For every $\rho\in[\rho_1,\rho_2]$, there exists a non-empty Aubry-Mather set of rotation number $\rho$.
Aubry-Mather sets defined as such can be obtained as limits of monotone Birkhoff periodic orbits \cite{Katok82}.

In \cite{Mather91}, Mather has proved that, given a bi-infinite sequence of action-minimizing Aubry-Mather sets inside a BZI, three exists a trajectory that visits arbitrarily close each of the Aubry-Mather set in the sequence, in the prescribed order. There is also a  topological version of this result,  due to Hall \cite{Hall1989}.

\subsection{Main result}
\label{subsection:main}

We consider a  $C^r$-differentiable Hamiltonian system, with $r$ sufficiently large,  given by a smooth function $H_\eps:\mathbb{R}^6\to\mathbb{R}$ of the form
\begin{equation}
H_\eps(p,q)=H_0(p,q)+\eps H_1(p,q),
\end{equation}
where $p=(p_1,p_2,p_3)$, $q=(q_1,q_2,q_3)$, and the unperturbed Hamiltonian $H_0$ is of the  form
\begin{equation}
H_0(p,q)=H_{00}(p_1,q_1,p_2,q_2)+H_{01}(p_3,q_3).
\end{equation}
On $\mathbb{R}^6$ we consider the canonical symplectic form $\omega=\sum_{i=1}^{3}dq_i\wedge dp_i$. Denote by $\textbf{0}_4=(0,0,0,0)$ the origin of the $(p_1,q_1,p_2,q_2)$-phase-space, and by $\textbf{0}_2=(0,0)$ the origin of the $(p_3,q_3)$-phase-space.

Suppose that $H_{00}(\textbf{0}_4)=H_{01}(\textbf{0}_2)=0$, and that there exists $c_1>0$, such that each $c\in(0,c_1]$   is a regular value for $H_{00}$. Let \[\Lambda_{0,c}=\{(p_1,q_1,p_2,q_2)\,|\,H_{00}(p_1,q_1,p_2,q_2)=c\}\] be   the corresponding energy hypersurface  of $H_{00}$ in the $(p_1,q_1,p_2,q_2)$-phase-space.

We now make some assumptions on the Hamiltonian $H_\eps$ which will be used in formulating the main result.

\begin{itemize}
\item [(A1)] The Hamiltonian   $H_{00}$ satisfies the following properties:
\begin{itemize}
\item[(i)]   $H_{00}$ has a non-degenerate minimum point at the origin $\textbf{0}_4$.

\item [(ii)]  The following strict convexity condition holds: there exists $C>0$ such that $(D^2H_{00}-C\cdot\textrm{id})$ is positive definite at all points $(p_1,q_1,p_2,q_2)\in\mathbb{R}^4$.

\item[(iii)]  $H_{00}$ is  Liouville integrable:  there exists a first integral $K$ such that $\{H_{00},K\}$ are functionally independent almost everywhere and in involution.

\item [(iv)] $K$ satisfies a  Bott condition on $\Lambda_{0,c}$, and   the  critical submanifolds of $K$ on $\Lambda_{0,c}$  are only circles of either elliptic or hyperbolic type, for each $c\in(0,c_1]$.

\item [(v)]  On each maximal action-angle domain $\Omega_0$ with corresponding action-angle coordinate system $(I_1,I_2,\phi_1,\phi_2)\in A\times\mathbb{T}^2$,   the Hamiltonian $H_{00}$ satisfies the following condition:  $(I\mapsto H_{00}(I))^{-1}(c)$  is a strictly-convex hypersurface  in $A$,
    for each $c\in (0,c_1]$.

\end{itemize}
\item [(A2)]  The Hamiltonian $H_{01}$  satisfies the following properties:
\begin{itemize}
\item[(i)]  $H_{01}$  has  a non-degenerate saddle point at $\textbf{0}_2$, and the set \[\gamma=\{(p_3,q_3)\,|\, H_{01}(p_3,q_3)=0\}\] is the union of two smooth curves $\gamma^-,\gamma^+$ (separatrices) that meet only at $\textbf{0}_2$, and contain no other critical point of $H_{01}$. Each $\gamma^\pm$ represents the stable and unstable manifold of $\textbf{0}_2$; the two manifolds coincide.
\item [(ii)]  The absolute values of the Lyapunov exponents of the flow of $H_{00}$ on each $\Lambda_{0,c}$, $c\in(0,c_1]$ are less than $\mu$, where $\pm\mu$ denote the Lyapunov exponents of the saddle point $\textbf{0}_2$ for the flow of $H_{01}$.

\end{itemize}
\end{itemize}

Now we discuss the above assumptions.

Assumptions  (A1-i) and   (A1-ii) imply that $\textbf{0}_4$ is the unique critical point of $H_{00}$ within some neighborhood of $\textbf{0}_4$, and that $H_{00}$ is a Morse function in some neighborhood of $\textbf{0}_4$ in the $(p_1,q_1,p_2,q_2)$ phase-space. The Morse index of $\textbf{0}_4$ equals $0$. 
By the Morse Lemma, it follows that all $H_{00}$ level sets $\Lambda_{0,c}=\{H_{00}=c\}$, $c\in(0,c_1]$, are three-spheres in the $(p_1,q_1,p_2,q_2)$-phase-space, provided $c_1$ is sufficiently small.

Assumption (A1-ii) further implies that each energy hypersurface $\Lambda_{0,c} \simeq S^3$ in $\mathbb{R}^4$  is star-shaped relative to the origin, so the Hamiltonian flow on it is equivalent to the Reeb flow on a tight three-sphere. By Theorem \ref{thm:HWZ98} there exists a disk-like global surface of section $\mathcal{D}_{0,c}$ for the Hamiltonian flow restricted to $\Lambda_{0,c}$, and a corresponding return map $f_{0,c}$ that is conjugate to an area preserving map. See Appendix~\ref{subsection:global}.

Assumption  (A1-iii)  says that $H_{00}$ is  Liouville integrable, hence the complement of the critical set of the momentum map $\{H_{00},K\}$ is an open dense set that is a countable union of action-angle domains $\Omega_0^j$, $j\in\mathbb{N}$.

Assumption (A1-iv) says that for each $c\in(0,c_1]$, the first integral $K$ satisfies a Bott condition on each $\Lambda_{0,c}$, hence the complement of the critical set of $K$ on $\Lambda_{0,c}$ is in fact a finite union of open sets on which $H_{00}$ and $K$ are functionally independent. This implies that only finitely many of the action-angle domains from above, say $\Omega_0^j$, $j=1,\ldots,k$,  intersect  each $\Lambda_{0,c}$.
The sets $\Omega_{0,c}^j=\Omega^j_0\cap \Lambda_{0,c}$, $j=1,\ldots,k$,  are mutually disjoint open sets, and  the complement of $\bigcup_{j=1,\ldots,k}\Omega^j_{0,c}$ is a nowhere dense set in $\Lambda_{0,c}$. Each  domain $\Omega_{0,c}^j$ can be described by a system $(I_1,\phi_1,\phi_2)$ of one action and two angle coordinates, induced by the corresponding action-angle coordinate system from (A1-v);  the second action coordinate $I_2$ is implicitly defined by the energy condition.
Thus, $\Lambda_{0,c}$ is the union of the critical sets of $K$ on $\Lambda_{0,c}$, which, by assumption, are only elliptic circles    or hyperbolic circles together with their stable and unstable manifolds, and of the domains $\Omega_{0,c}^j$, $j=1,\ldots,k$.
Each such domain $ \Omega_{0,c}^j$, $j=1,\ldots,k$ determines a corresponding domain on the global surface of section $\mathcal{D}_{0,c}$ for the Hamiltonian flow restricted to $\Lambda_{0,c}$, denoted  by
  $\widehat\Omega_{0,c}^j= \Omega_{0,c}^j\cap \mathcal{D}_{0,c}$. We can parametrize each
$\widehat\Omega_{0,c}^j$ by an action-angle coordinate system $(J,\theta)$, where $J$ can be taken equal to $I_1$, and $\theta$ is the angle coordinate symplectically conjugate to $J$ relative to the symplectic form induced on
$\widehat\Omega_{0,c}^j$.

Assumption (A1-v) is equivalent in two-degrees of freedom with the following iso-energetic non-degeneracy condition
\[ \det \begin{pmatrix}
           \displaystyle\frac{\partial \omega}{\partial I} & \omega \\[7pt]
           \omega^T & 0 \\
         \end{pmatrix}
       \neq 0,\]
for all $I\in A$, where ${(\cdot)}^T$ denotes the transpose of a matrix. See \cite{Lochak1992}.
This condition  can be also reformulated in terms of the frequency
map  $I\mapsto W(I)=[\omega_1(I):\omega_2(I)]\in P^1(\mathbb{R})$, where $\omega_1=\partial H /\partial I_1$, $\omega_2=\partial H /\partial I_2$. For each $I\in A$, $W(I)$ represents
the frequency ratio on the torus determined by $I$. The
iso-energetic non-degeneracy condition is equivalent to the condition
$\partial W/\partial I\neq 0$ for all $I\in A$, where the derivative is taken at constant energy; hence the frequency ratio map $I\mapsto W(I)$ from the energy hypersurface to the projective line $P^1(\mathbb{R})$ has maximal rank. The iso-energetic non-degeneracy condition implies that the rotation number changes between neighboring tori.  This condition can also be reformulated in terms of transversality, see \cite{DelshamsG96}.
Conditions (A1-v) combined with (A1-ii)  imply that the return map  $f_{0,c}$ to the global surface of section $\mathcal{D}_{0,c}$   is a twist map relative to the  action-angle coordinate system $(J,\theta)$, on each  $\widehat\Omega_{0,c}^j$, $j=1,\ldots,k$.

Assumption (A2) implies that  
\begin{equation}\label{eqn:nhim_0} \Lambda_0=\bigcup_{c\in[0,c_1]}\Lambda_{0,c}\subseteq \{(p,q)\in\mathbb{R}^6\,|\,p_3=q_3=0\}\end{equation} is a $4$-dimensional normally hyperbolic invariant manifold (with boundary) for the Hamiltonian flow. We also have that $W^u(\Lambda_0)=W^s(\Lambda_0)$, and that $\Lambda_0$ is foliated by the three-spheres $\Lambda_{0,c}$.
The next proposition says that, when the perturbation  applied to the system is small enough, the normally hyperbolic invariant manifold survives, and some part of it remains foliated by three-spheres.

\begin{prop}\label{prop:nhim_eps}
There exist $0<\eps_1$ and $0<c_0<c_1$ sufficiently small, such that    for each  $0<\eps<\eps_1$ there exists a normally hyperbolic locally invariant manifold $\Lambda_\eps$ that is $\eps$-close to $\Lambda_0$ in the $C^{r-2}$ topology, and  for each    $c\in[c_0, c_1]$,  $\Lambda_{\eps,c}:=\Lambda_\eps\cap \{H_\eps=c\}$ is a normally hyperbolic  invariant manifold diffeomorphic to a three-sphere. Moreover, there exists a parametrization $k_\eps:\Lambda_0\to\Lambda_\eps$ such that $k_\eps(\Lambda_{0,c})=\Lambda_{\eps,c}$ for all $0<\eps<\eps_1$ and $c\in[c_0,c_1]$.
\end{prop}

The proof of this proposition is given in Section \ref{section:proofs}.

Let $\Omega_\eps^j=k_\eps(\Omega_0^j)$, $j=1,\ldots,k$. Using the parametrization $k_\eps$ we can carry the coordinate system from  each domain $\Omega_0^j$ to a coordinate system on  $\Omega_\eps^j$, $j=1,\ldots,k$, which, with an abuse of notation, we still denote by $(I_1,I_2,\phi_1,\phi_2)$.
However, since the perturbed system is not necessarily integrable, the action level sets are not in general preserved by the Hamiltonian flow of $H_\eps$.
Also note that the boundaries of the  domains $\Omega^j_\eps$  are no longer invariant sets, as they may get destroyed by the perturbation.
Similarly, each domain $\Omega^j_{0,c}$ in $\Lambda_{0,c}$ corresponds via the restriction of $k_{\eps}$ to $\Lambda_{0,c}$ to a domain $\Omega^j_{\eps,c}$ in $\Lambda_{\eps,c}$. Thus, there exists a  coordinate system  $(I_1,\phi_1,\phi_2)$ induced via $k_\eps$ on $\Omega^j_{0,c}$ by the corresponding coordinate system on 
$\Omega^j_{0,c}$.

For a perturbation $H_1$ of a `generic type', and for all $\eps$ sufficiently small, $W^u(\Lambda_\eps)$ and $W^s(\Lambda_\eps)$ intersect transversally. More precisely, we will show that there exists an open and dense set of perturbations $H_1$ for which there exists a collection of homoclinic channels $\Gamma_{\eps}^{i}$, $i=1,\ldots, n$, and corresponding scattering maps $S_{\eps}^{i}:U^i_{\eps}\to V^i_{\eps}$ associated to $\Gamma_{\eps}^{i}$, with $U^i_\eps \subseteq \Lambda_\eps$,  such that,
for every $j\in\{1,\ldots k\}$, each level set of $I_1$ in $\Omega^j_\eps$ intersects some $U^i_\eps$.
The notion of a scattering map is recalled in Appendix~\ref{app:scattering}. The above condition on the homoclinic channels is spelled out in condition (A3), Section \ref{sec:geompert}. The genericity of this condition is shown in Section \ref{subsub:generic}.

In order to formulate the main result of this paper, we need to discuss the existence of a disk-like global surfaces of section  for the flow of $H_\eps$ restricted to $\Lambda_{\eps,c}$, $c\in[c_0, c_1]$.

\begin{prop}\label{prop:globals} There exists $\eps_1>0$ such that for each $0<\eps<\eps_1$ and each $c\in[c_0, c_1]$ there exists a family of disk-like global surfaces of sections $\mathcal{D}_{\eps,c}$ for the Hamiltonian flow of $H_\eps$ restricted to $\Lambda_{\eps,c}\simeq S^3$, which depends smoothly  on $c$ and $\eps$.
The return map $f_{\eps,c}$ to $\mathcal{D}_{\eps,c}$ is a twist map and is conjugate to an area preserving map.
\end{prop}

The proof of this proposition is given in Section \ref{section:proofs}.

Each  domain $\Omega^j_{\eps,c}$,  $j=1,\ldots,k$, determines an open domain $\widehat\Omega^j_{\eps}=\Omega^j_{\eps,c}\cap \mathcal{D}_{\eps,c}$   in $\mathcal{D}_{\eps,c}$. Each  $\widehat\Omega^j_{\eps,c}$ is diffeomorphic to an annulus (punctured disk) and can be described in terms of a system of action-angle coordinates $(J,\theta)$   induced by the corresponding system on $\widehat\Omega_{0,c}^j$. We use the same notation for the induced coordinate system as for the original one. The twist condition in Proposition \ref{prop:globals} is understood relative to the induced action-angle coordinates $(J,\theta)$ on
$\widehat\Omega^j_{\eps,c}$.

\begin{prop}\label{prop:kam}
There exists $\eps_1>0$ sufficiently small, such that, for each $\eps\in(0,\eps_1)$ and each $c\in[c_0,c_1]$, on each domain $\widehat{\Omega}^j_{\eps,c}$, $j=1,\ldots, k$, the map $f_{\eps,c}$   has a Cantor set of KAM invariant tori $\{\mathcal{T}_\alpha\}_{\alpha\in\mathcal{A}_j}$   that survive from  the  foliation by  $f_{0,c}$-invariant tori  of the domain $\widehat{\Omega}^j_{0,c}$.
\end{prop}

The proof of this proposition is given in Section \ref{section:proofs}. Besides KAM tori, since $f_{\eps,c}$ is a monotone twist map, and is conjugate to an area preserving map, for every region in  $\widehat{\Omega}^j_{0,c}$ bounded by a pair of invariant KAM tori $\mathcal{T}_{1}$, $\mathcal{T}_{2}$, and for each rotation number $\rho$ in the rotation interval $[\rho_1,\rho_2]$ determined by the boundary tori, there exists an Aubry-Mather set of that rotation number. See Section \ref{subsection:aubrymather}.

The main theorem below aims to describe  the existence of `global diffusion' relative to the three-sphere $\Lambda_{\eps,c}$.
We give a precise sense of  `global diffusion' as follows.  Choose   a global section $\mathcal{D}_{\eps,c}$ of the Hamiltonian flow restricted to $\Lambda_{\eps,c}$. Inside  $\mathcal{D}_{\eps,c}$, choose a collection of invariant subsets $\{\mathcal{R}_1, \mathcal{R}_2,\ldots,\mathcal{R}_n\}$ so that each  $\mathcal{R}_i$ is either a KAM invariant torus or an Aubry-Mather set within  some domain $\widehat\Omega^j_{\eps,c}$, for $j\in\{1,\ldots,k\}$.

\begin{thm}\label{thm:main} Assume that $H_\eps$ satisfies the assumptions (A1)-(A2) from above, and assumption (A3) from Section \ref{sec:geompert}.  Then, there exist  an open and dense set of smooth perturbations $H_1$ and  $\eps_1>0$  sufficiently small, such that, for every $\eps\in(0,\eps_1)$, every $c\in[c_0, c_1]$,   given a global section $\mathcal{D}_{\eps,c}$ of the Hamiltonian  flow of $H_\eps$ restricted to $\Lambda_{\eps,c}$ and  a collection of subsets $\{\mathcal{R}_1, \mathcal{R}_2,\ldots,\mathcal{R}_n\}$  of $\mathcal{D}_{\eps,c}$ as above,  and a positive real $\delta>0$, there exists a trajectory $\phi^t_\eps(x)$  of the Hamiltonian flow and a sequence of times $t_1<t_2<\ldots<t_n$ such that $d(\phi^{t_i}_\eps(x),\mathcal{R}_i)< \delta$ for $i=1,\ldots,n$.
\end{thm}

The proof of Theorem \ref{thm:main} is given in Section \ref{section:proofs}.

\begin{rem} From the statement of Theorem \ref{thm:main} it may appear that the choice of the invariant objects $\{\mathcal{R}_1, \mathcal{R}_2,\ldots,\mathcal{R}_n\}$ depends on the choice of the global section $\mathcal{D}_{\eps,c}$. This is not true.  A set $\mathcal{R}_i$ that is a $1$-dimensional KAM torus for the map $f_{\eps,c}$ in
$\mathcal{D}_{\eps,c}$ gives rise to a  $2$-dimensional KAM torus    for the flow $\phi^t_\eps$. Also,  a set $\mathcal{R}_i$ that is an Aubry-Mather sets for the map $f_{\eps,c}$ in
$\mathcal{D}_{\eps,c}$ gives rise to a lamination over the Aubry-Mather set for the the flow $\phi^t_\eps$. (An interesting analysis of these laminations appear in \cite{Fayad2008}).
Thus, the above theorem yields trajectories that shadow a prescribed collection of invariant objects that is independent of the surface of section.

\end{rem}

\begin{rem}
The assumption in (A1-iv) that the critical submanifolds of $K$ on $\Lambda_{0,c}$ are only circles does not seem to be essential; we make it mainly to simplify the geometric picture. If a critical submanifold is a torus, then there is a neighborhood of the torus  which is an action-angle domain. Also, if a critical submanifold is a  Klein bottle, then there is
neighborhood of the Klein bottle where the dynamics can be reduced, via a two-sheet covering, to the dynamics near a critical torus. See   \cite{Marco12}.
\end{rem}

\begin{rem}
For most of the argument, the differentiability class $C^r$ of the Hamiltonian can be chosen with $r$  fairly low, e.g. $r\geq 5$. The existence of a global surface of sections for Hamiltonian flows on   three-dimensional strictly convex
energy surfaces is proved in \cite{HZW98} for $C^\infty$-differentiable systems. In a personal communication, H. Hofer informed us that the result  remains valid for  $C^r$-differentiable systems, with $r$ sufficiently large.
\end{rem}

\begin{rem}
The strict convexity condition (A1-ii) is not a symplectic condition. However it can be replaced by a symplectic condition in terms of the contact structure induced on each three-sphere $\Lambda_{0,c}$.   More precisely, one can replace it by the following conditions: (i)~the three-sphere is of contact type, (ii) the contact structure is tight, (iii) the Reeb flow associated to the contact structure is equivalent to the Hamiltonian flow restricted to the three-sphere, (iv) every periodic orbit of the Reeb flow has Conley-Zehnder index greater than or equal to three. Under these conditions, a deep result from \cite{HZW98} says that the Reeb flow has a disk-like global surface of section, on which the Poincar\'e return map is equivalent to an area preserving map. See Appendix~\ref{subsection:global}.
\end{rem}

\section{Examples}\label{sec:application}
We consider a Hamiltonian system consisting of the coupling of two non-harmonic oscillators and one pendulum subject to a small perturbation, given by \[H_\eps(p,q)=H_{00}(p_1,q_1,p_2,q_2)+H_{01}(p_3,q_3)+\eps H_1(p,q),\] where
\begin{equation}\begin{split}
H_{00}(p_1,q_1,p_2,q_2)&= a_1\frac{p_1^2+q_1^2}{2} +a_2\frac{p_2^2+q_2^2}{2} + \frac{b_1}{2} \left(\frac{p_1^2+q_1^2}{2}\right)^2+\frac{b_2}{2} \left(\frac{p_2^2+q_2^2}{2}\right)^2, \\
H_{01}(p_3,q_3)&= \frac{1}{2}p_3^2+\lambda^2(\cos q_3-1) ,
\end{split}
\end{equation}with $a_1,a_2,b_1,b_2>0$, $a_1\neq a_2$, $b_1\neq b_2$, $\lambda>0$.

We note that $H_0=H_{00}+H_{01}$ is an integrable Hamiltonian. 

The   Hamilton equations corresponding to $H_\eps$ are:
\begin{equation}\label{eqn:hamex}\begin{split}
\dot q_1=&a_1p_1+b_1p_1\frac{p_1^2+ q_1^2}{2}+\eps \frac{\partial H_1}{\partial p_1},\\
\dot p_1=&-a_1q_1-b_1q_1\frac{p_1^2+ q_1^2}{2}-\eps\frac{\partial H_1}{\partial q_1},\\
\dot q_2=&a_2p_2+b_2p_2\frac{p_2^2+ q_2^2}{2}+\eps \frac{\partial H_1}{\partial p_2},\\
\dot p_2=&-a_2q_2-b_2q_2\frac{p_2^2+ q_2^2}{2}-\eps\frac{\partial H_1}{\partial q_2},\\
\dot q_3=&p_3+\eps \frac{\partial H_1}{\partial p_3},\\
\dot p_3=&\lambda^2 \sin q_3-\eps\frac{\partial H_1}{\partial q_3},
\end{split}\end{equation}

First, we let $\eps=0$. The point $(p,q)={\bf 0}_{6}$ is an equilibrium point of the Hamiltonian flow, of linearized type center$\times$center$\times$saddle.
The eigenvalues of the linearized system at this point are $\pm a_1 i,\pm a_2 i,\pm\lambda$.
The Hamiltonian $H_{01}$ satisfies the conditions on (A2).
The Hamiltonian $H_{00}$ has a non-degenerate minimum at ${\bf 0}_4$, as in (A1-i). 

We now check (A1-ii).
We compute
\[D^2H_{00}=\left(\begin{array}{ll }D_1 &0 \\
0 & D_2
\end{array}
\right),\]
where \[D_i=\left(\begin{array}{ll}a_i+\frac{b_i}{2}(3p_i^2+q_i^2) & b_ip_iq_i \\
b_ip_iq_i &a_i+\frac{b_i}{2}(p_i^2+3q_i^2)\end{array}
\right),\]
with $i=1,2$.

We have \begin{equation*}
\begin{split}\langle D^2H_{00}&(u_1,v_1,u_2,v_2),(u_1,v_1,u_2,v_2)\rangle=\\
&=\sum_{i=1,2} a_i(u_i^2+v_i^2)+\frac{b_i}{2}[ (3p_i^2+q_i^2)u_i^2+4p_iq_iu_iv_i+ (p_i^2+3q_i^2)v_i^2]\\
&\geq \sum_{i=1,2} a_iu_i^2+a_iv_i^2\\
&\geq C\sum_{i=1,2}   u_i^2+  v_i^2\\
&=C\langle  (u_1,v_1,u_2,v_2),(u_1,v_1,u_2,v_2)\rangle,
\end{split}
\end{equation*}
provided we choose $0<C<\min\{a_1,a_2\}$.

It is clear that the Hamiltonian $H_{00}$ is Liouville integrable integrable, with two first integrals $K_1(p,q)= \frac {p_1^2+q_1^2}{2}$ and $K_2(p,q)=\frac{p_2^2+q_2^2}{2}$, thus verifying (A1-iii). Of course $H_{00}=h_1(K_1)+h_2(K_2)$, where $h_i(z)=a_iz+\frac{b_i}{2}z^2$, $i=1,2$, is also a first integral.
The value $c=0$ is the only critical value of $H_{00}$.
If $c>0$, the set $\Lambda_{0,c}=\{(p,q)\in H^{-1}_0(c)\,|\, p_3=q_3=0\}$ is invariant under the Hamiltonian flow of $H_0$, and $\Lambda_{0,c}\simeq S^3$.

We fix $c>0$ sufficiently small.
The critical submanifolds of $\Lambda_{0,c}$ are the circle $\chi^1_{0,c}$ given by $p_1^2+q_1^2=\frac{-a_1+\sqrt{a_1^2+2b_1c}}{b_1}$, $p_2=q_2=0$, and the circle $\chi^2_{0,c}$ $p_2^2+q_2^2=\frac{-a_2+\sqrt{a_2^2+2b_2c}}{b_2}$, $p_1=q_1=0$. Both circles represent closed orbits of elliptic type, thus verifying (A1-iv). Note that the two circles form a Hopf link.

We make a canonical change of coordinates to  action-angle  coordinates
\begin{eqnarray*}I_1=\frac{1}{2}(q_1^2+p_1^2),\,\phi_1=\tan^{-1} (p_1/q_1),\\
I_2=\frac{1}{2}(q_2^2+p_2 ^2),\,\phi_2=\tan^{-1} (p_2/q_2),
\end{eqnarray*}
so that $\sum_{i=1,2}{dq_i\wedge dp_i}=\sum_{i=1,2}dI_i\wedge d\phi_i$.

The Hamiltonian $H_{00}$ in action-angle coordinates is
\[H_{00}(I_1,\phi_1,I_2,\phi_2)=a_1I_1+a_2I_2+\frac{b_1}{2}I_1^2+\frac{b_2}{2}I_2^2.\]
The  equations of motion on $\Lambda_0$ in action-angle coordinates are
\begin{equation}\label{eqn:hamex0}\begin{split}
\dot I_1=&0,\\
\dot \phi_1=&a_1+b_1I_1,\\
\dot I_2=&0,\\
\dot \phi_2=&a_2+b_2I_2.
\end{split}
\end{equation}
We have
\[\det\left(
        \begin{array}{cc}
          \frac{\partial^2H_{00}}{\partial^2I} & \frac{\partial H_{00}}{\partial I}  \\
          \left(\frac{\partial H_{00}}{\partial I}\right)^T& 0 \\
        \end{array}
      \right)=-b_2(a_1+b_1I_1)^2-b_1(a_2+b_2I_2)^2.\]
Since $a_1,a_2,b_1,b_2>0$, this determinant never vanishes, so the iso-energetic non-degeneracy condition is satisfied, thus verifying (A1-v).

\begin{figure}
\centering
\includegraphics[width=0.25\textwidth]{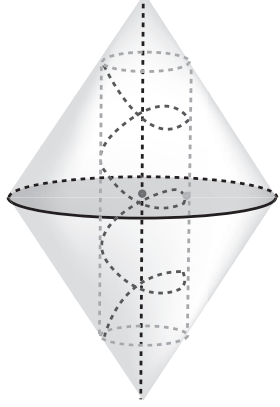}
\caption{Disk-like global surface of section for the Hamiltonian flow on the three-sphere.}
\label{fig:threesphere0}
\end{figure}

Due to (A1-ii),   $\Lambda_{0,c}$ bounds a strictly convex domain, and  Theorem \ref{thm:HWZ98} tells us that there exists a disk-like global surface of section $\mathcal{D}_{0,c}$ of the flow of $H_{00}$. However, in this example  such a global surface of section is  given explicitly  by $p_2=q_2=0$, which is the disk bounded by the closed orbit given by  $p_1^2+q_1^2=\frac{-a_1+\sqrt{a_1^2+2b_1c}}{b_1}$, $p_2=q_2=0$.  See Fig. \ref {fig:threesphere0}. This disk is foliated by invariant circles of the type $\{I_1=\textrm{const}\}$, with $I_2$ implicitly given by the energy condition $H_{00}(I_1,I_2)=c$.

From the above assumptions, it follows that the statements from Proposition \ref{prop:nhim_eps}, \ref{prop:globals}, and \ref{prop:kam} hold true. In particular, there exist $0<c_0<c_1$ and $0<\eps_1$ such that, for each $c\in[c_0,c_1]$ and each $0<\eps<\eps_1$, the energy hypersurface $\{H_\eps=c\}$ contains a normally hyperbolic invariant manifold $\Lambda_{\eps,c}\simeq S^3$. This can be endowed with  a coordinate system $(I_1,\phi_1,\phi_2)$ induced by the corresponding coordinate system on $\Lambda_{0,c}$.  
There exists a 
disk-like global surface of section $\mathcal{D}_{\eps,c}\subseteq\Lambda_{\eps,c}$. The boundary $\partial \mathcal {D}_{\eps,c}$ is a periodic orbit $\chi^1_{\eps,c}$. There exists another periodic orbit $\chi^2_{\eps,c}$ in ${\Lambda}_{\eps,c}$ that is transverse to $\textrm{int}(\mathcal{D}_{\eps,c})$, forming a Hopf link with $\chi^1_{\eps,c}$.
Let $p^2_{\eps,c}$ be the intersection point between $\chi^2_{\eps,c}$ with $\mathcal{D}_{\eps,c}$.   Let $f_{\eps,c}:\textrm{int}(\mathcal {D}_{\eps,c})\to \textrm{int}(\mathcal {D}_{\eps,c})$ be the Poincar\'e return map to $\mathcal {D}_{\eps,c}$  for the flow restricted to $\Lambda_{\eps,c}$. The return map restricted to the open annulus $A_{\eps,c}=\mathcal {D}_{\eps,c}\setminus\{p^2_{\eps,c}\}$  is equivalent to an area preserving map. It also satisfies a monotone twist condition. Hence we can select in $A_{\eps,c}$ some finite, arbitrary collection of invariant sets $\{\mathcal{R}_1,\ldots,\mathcal{R}_n\}$, such that each $\mathcal{R}_i$ is either a KAM invariant torus or an Aubry-Mather set for $f_{\eps,c}$.

Assume  that $H_1$ satisfies  the  non-degeneracy condition (A3) from Section \ref{sec:geompert}. 
Then Theorem \ref{thm:main} applies. 
Under these assumptions, the stable and unstable manifold $W^u(\Lambda_{\eps,c})$ and $W^u(\Lambda_{\eps,c})$ intersect transversally, and there exists a finite collection of homoclinic channels $\{\Gamma_{\eps,c}^i\}$, and corresponding scattering maps $S_{\eps,c}^i:U_{\eps,c}^i\to V_{\eps,c}^i$, $i=1,\ldots,n$, such that $\bigcup_{i=1,\ldots,n} U_{\eps,c}^i$ covers the whole range of the action variable $I_1$ on $\Lambda_{\eps,c}$. That is, $\Pi_{I_1}(\bigcup_{i=1,\ldots,n} U_{\eps,c}^i)\supseteq \Pi_{I_1}(\Lambda_{\eps,c}\setminus (\chi^1_{\eps,c}\cup \chi^2_{\eps,c}))$.

Under these conditions, Theorem \ref{thm:main} implies that the Hamiltonian flow  has orbits that visits arbitrarily closely the prescribed collection of invariant subsets    $\{\mathcal{R}_1,\ldots,\mathcal{R}_n\}$ in $\mathcal {D}_{\eps,c}$.

\begin{rem}\label{rem:RTBP}
We compare this example with the  spatial circular restricted three-body problem. In a co-rotating coordinate system, with the primaries of masses $1-\mu$, $\mu$ placed at $-\mu$, $1-\mu$, respectively, the motion of the infinitesimal mass is described by the following Hamiltonian
\[H(x,y,z,\dot x, \dot y,\dot
z)=\frac{1}{2}(\dot x ^2 +\dot y^2+\dot z^2)-\omega (x,y,z),\]
where the `effective potential' $\omega$ is given by
\[\omega
(x,y,z)=\frac{1}{2}(x^2+y^2)+ \frac{1-\mu}{r_1}+\frac{\mu}{r_2},\]
with $r_1=((x+\mu)^2+y^2+z^2)^{1/2}$ and $r_2=((x-1+\mu)^2+y^2+z^2)^{1/2}$.
This problem has five equilibria,   $L_1,L_2,L_3$ of  (center$\times$center$\times$ saddle)-type, and
$L_4,L_5$ of   (center$\times$center$\times$ center)-type.  We denote by $L_1$  the equilibrium between the primaries. We will assume that $\mu$ is very small.

After performing a translation of coordinates and using  a normal form about $L_1$, the Hamiltonian can be written in a neighborhood of $L_1$ in the form
\[H (p,q)=H_{00}(p_1,q_1,p_2,q_2)+H_{01}(p_3,q_3)+H_1(p,q),\] with $H_{00},H_{01},H_1$ of the forms
\begin{equation*}
\begin{split}
H_{00}(p_1,q_1,p_2,q_2)&=\frac{1}{2}\omega_1(p_1^2+q_1^2)+\frac{1}{2}\omega_2(p_2^2+q_2^2),\\
H_{01}(p_3,q_3)&=\lambda p_3q_3,\\
H_1(p_1,q_1,p_2,q_2,p_3,q_3)&=\sum_{k\geq 3}H_1^k(p_1,q_1,p_2,q_2,p_3,q_3),
\end{split}
\end{equation*}
where $\lambda, \omega_1,\omega_2$ are positive reals, and $H_1^k(p,q)$, $k\geq 3$,  is a  homogeneous
polynomial of order $k$ in the variables $p_i,q_i$, $i = 1, \ldots, 3$.

There exists a center manifold $W^c$ that is tangent at $L_1$ to the subspace corresponding to the imaginary eigenvalues $\pm i\omega_1,\pm i\omega_2$. Every energy level of $H$ sufficiently close to $L_1$ intersects $W^c$ along a three-sphere $\Lambda_c\simeq S^3$. See, e.g., \cite{Moeckel2004}.
It is easy to see  that $H_{00}$ satisfies the  strict convexity condition in a neighborhood of $L_1$. 
In \cite{Sama04} it was proved that the iso-energetic non-degeneracy condition is satisfied on $\Lambda_c$, provided $c$ is sufficiently close to the energy level of $L_1$. 
Again, we do not have to invoke Theorem \ref{thm:HWZ98} to prove the existence of a disk-like global surface of section, since the  $z$-axis symmetry of the problem ensures that
the hyperplane $\{z=0\}$ is a global surface of section for the Hamiltonian flow restricted to~$\Lambda_c$.

A  significant  difference from the previous example is that the term $H_1$ is intrinsic to the problem, and cannot be regarded  as a small, `generic' type of perturbation. Thus, we cannot apply Theorem \ref{thm:main}. This is a typical
shortcoming of `generic' type of diffusion mechanisms, namely that in many instance such mechanisms cannot be applied to concrete examples. (It is well known that there are examples of `generic'  systems for which not a single example  is known.)

In order to establish the existence of orbits that exhibit global diffusion relative to $\Lambda_c$  in this problem  one needs   to explicitly find (numerically or analytically) a collection of homoclinic intersections and corresponding scattering maps that can be used to achieve  global diffusion. Some numerical evidence is given  in \cite{DelshamsGR10,DelshamsGR13,Sama04}.
\end{rem}

\section{Proofs of the results}\label{section:proofs}
\subsection{Proof of Proposition \ref{prop:nhim_eps}}
As we already mentioned in Section \ref{subsection:main}, by  (A1-i) and   (A1-ii),  $\textbf{0}_4$ is the unique critical point of $H_{00}$ within some neighborhood of $\textbf{0}_4$,  $H_{00}$ is a Morse function in some neighborhood of $\textbf{0}_4$ in the $(p_1,q_1,p_2,q_2)$ phase-space, and  the  Morse Lemma implies that the level sets $\Lambda_{0,c}=\{H_{00}(p_1,q_1,p_2,q_2)=0\}$ 
are three-spheres in the $(p_1,q_1,p_2,q_2)$-phase-space, provided $c_1$ is sufficiently small.

By (A2) we obtain  that $\Lambda_0=\bigcup_{c\in[0,c_1]}\Lambda^c_0$  is a compact, normally hyperbolic invariant manifold with boundary for the Hamiltonian flow of $H_0$. The persistence theorem  of normally hyperbolic invariant manifolds  implies that, if the perturbation applied to $H_0$ is sufficiently
small, then $\Lambda_0$ is survived by a normally hyperbolic locally invariant manifold $\Lambda_\eps$ for the Hamiltonian flow of $H_\eps$. Recall that local invariance means that there exists a neighborhood $\mathscr{U}_{\Lambda_\eps}$ of $\Lambda_\eps$ such that any trajectory of the Hamiltonian flow that stays in $\mathscr{U}_{\Lambda_\eps}$ is contained in $\Lambda_\eps$; equivalently, the Hamiltonian vector field of $H_\eps$ is tangent to $\Lambda_\eps$ at all   points. Moreover, there exists a smooth parametrization $k_\eps:\Lambda_0\to\Lambda_\eps$ of $\Lambda_\eps$ that depends smoothly on $\eps$, and with $k_0=\textrm{id}_{\Lambda}$. The manifold $\Lambda_\eps$ is $\eps$-close to $\Lambda_0$ in the $C^{r-2}$-topology.
The parametrization $k_\eps$ is not unique; if we compose $k_\eps$ to any diffeomorphism of $\Lambda_0$ we obtain a new parametrization. See  \cite{HirschPS77,DelshamsLS2006}. 

We consider the restriction of the function $H_\eps$ to $\Lambda_\eps$. 
By the implicit function theorem, for every sufficiently small $\eps>0$  there exists a critical point $\textbf{o}_\eps$ of $H_\eps$ on $\Lambda_\eps$ which is a nondegenerate minimum point  of $H_{\eps}$. Fixing $\eps_1>0$ sufficiently small we can choose a neighborhood $\mathscr{V}$ of $\textbf{o}_\eps$ in $\mathbb{R}^6$ such that $\textbf{o}_\eps$ is the unique critical point of $H_\eps$ on $\Lambda_\eps\cap\mathscr{V}$, for all $0\leq \eps<\eps_1$. 

For each $0\leq\eps<\eps_1$, there exist  $0<c_0<c_1$, such that the following hold true:
\begin{itemize}
\item for all $0\leq c<c_1$, the level sets $\{H_\eps=c\}$ in $\Lambda_\eps$ are contained in $\mathscr{V}$;
\item for each $0<c\leq c_0$, the Morse Lemma applies and all level sets $\{H_\eps=c\}$ in $\Lambda_\eps$ are  three-spheres;
\item for each $c\in[c_0,c_1]$, the norm of the gradient   of $H_\eps$ restricted to $\Lambda_\eps$ is bounded away from zero, and all corresponding level sets $\{H_\eps=c\}$ in $\Lambda_\eps$ are diffeomorphic to one another, hence they are three-spheres.
\end{itemize}
Thus $\Lambda_{\eps,c}\simeq S^3$ for every $0<\eps<\eps_1$ and $c\in[c_0,c_1]$. The region $\{c_0\leq H_\eps\leq c_1\}\cap \Lambda_\eps$ is diffeomorphic to the cylinder $S^3\times [c_0,c_1]$, and is foliated by $\Lambda_{\eps,c}\simeq S^3$.
The   parametrization $k_\eps$ of $\Lambda_\eps$ introduced above  can be chosen to map each three-sphere $\Lambda_{0,c}$ of the unperturbed system into a three-sphere $\Lambda_{\eps,c}$ of the perturbed system.
\qedsymbol

\subsection{Proof of Proposition \ref{prop:globals}}
By Proposition \ref{prop:nhim_eps}, there exists a parametrization $k_\eps:\Lambda_0\to\Lambda_\eps$, with $k_0=\textrm{id}$, that induces a parametrization $k_{\eps}:\Lambda_{0,c}\to\Lambda_{\eps,c}$ for each~$c$.
Choose a parametrization of $\Lambda_{0,c}$ given by $z\in S^3\mapsto \kappa(z)= \sqrt{g(z)}z\in\Lambda_{0,c}$, with $g:S^3\to\mathbb{R}$ smooth. Let $\tilde k_\eps=k_\eps\circ \kappa$. Let $\lambda_\eps$ be the contact form on
$\Lambda_\eps$ and let $\tilde k_\eps^*\lambda_\eps$ be the contact form on $S^3$ given by the pull back of $\lambda_\eps$.
Since $\lambda$ is tight, we have that  $\tilde k_\eps^*\lambda_\eps$ is also tight provided $\eps$ is sufficiently small.

At this point on $S^3$ we have the standard contact form  $\lambda$, and the family of tight contact forms
$\tilde k_\eps^*\lambda_\eps$, with $\eps>0$ small. By Gray's Stability Theorem, there exists a smooth family of positive functions $\Psi_\eps$ such that $\tilde k_\eps^*\lambda_\eps=\sqrt{\Psi_\eps}\lambda$. We have that $\Psi_\eps$ is $C^{r-2}$-close to the identity for $\eps$ sufficiently small, and so the Reeb flow of $\tilde k_\eps^*\lambda_\eps$ is $C^{r-2}$-close to the  Reeb flow of $\lambda$. Since $\lambda$ is strictly convex, it is dynamically convex.
It then follows that $\tilde k_\eps^*\lambda_\eps$ is dynamically convex for all $\eps$ sufficiently small.

By virtue of Theorem \ref{thm:HWZ98}, we can choose $\eps_1$ sufficiently small such that, for all $0\leq\eps<\eps_1$, there exists a disk-like surface of section $\mathcal{D}_{\eps,c}$ for the flow of $H_\eps$ restricted to $\Lambda_{\eps,c}$. Moreover, since the property of being a global surface of section is an open condition, then we can choose the family of disks $\mathcal{D}_{\eps,c}$ to be $C^{r-2}$-smoothly depending on $\eps$ and $c$.

In addition, the boundary $\partial \mathcal {D}_{\eps,c}$ is a periodic orbit $\chi^1_{\eps,c}$. There exists another periodic orbit $\chi^2_{\eps,c}$ in ${\Lambda}_{\eps,c}$ transverse to $\textrm{int}(\mathcal{D}_{\eps,c})$, forming a Hopf link with $\chi^1_{\eps,c}$. See  Remark~\ref{rem:RTBP}.
\qedsymbol

\subsection{Proof of Proposition \ref{prop:kam}}
We start by noting that, in the unperturbed system, each domain $\widehat \Omega^{j}_{0,c}$ in $\mathcal{D}_{0,c}$, $j=1,\ldots,k$,   is an open annulus foliated by $1$-dimensional tori (essential circles) that are invariant under $f_{0,c}$, the first return map to  $\mathcal{D}_{0,c}$. Proposition \ref{prop:globals} says that $f_0$ is conjugate to an area preserving map, and  assumption (A1-v) that
$f_{0,c}$ is a twist map relative to the $(J,\theta)$ action-angle coordinate system on  $\widehat \Omega^{j}_{0,c}$. 
Provided $\eps_1$ is small enough, for each $0<\eps<\eps_1$,  $\widehat \Omega^{j}_{\eps,c}$ is an open annulus,     $f_\eps$ is conjugate to an area preserving map, and  is also a twist map relative to the action-angle coordinates  $(J,\theta)$ induced on  $\widehat \Omega^{j}_{\eps,c}$.
Applying Moser Twist Mapping Theorem \cite{Moser1973}\footnote{The Moser Twist Mapping Theorem was proved for $C^3$-mappings by M.R. Herman.}, we obtain that, for all $0<\eps<\eps_1$  with $\eps_1$ small enough, there exists a KAM family of tori invariant under $f_\eps$ that survives from the family of $f_{0,c}$-invariant tori that foliates $\widehat\Omega^j_{0,c}$.
\qedsymbol

\subsection{Proof of Theorem \ref{thm:main}}
\subsubsection{The geometry of the unperturbed system}\label{sec:geomunpert}
The assumption (A1-iv) is that all critical submanifolds of the first integral $K$ on $\Lambda_{0,c}$ are circles that are non-degenerate in the sense of Bott.  These circles are necessarily periodic orbits of the Hamiltonian flow. The  index of a  critical circle is the number of negative eigenvalues of the restriction of $D^2K$ on a subspace transverse to the circle.  These circles are also assumed to be of hyperbolic or elliptic type. In particular, the circles of index $0$ or $2$ are elliptic, and the circles of index $1$ are hyperbolic.

The connected components of the complement in $\Lambda_{0,c}$ of these  circles  are action-angle domains,
in the sense that on  each component $\Omega^j_{0,c}$ there exists a system of one action and two angle coordinates $(I_1,\phi_1,\phi_2)$.
More precisely,  for a critical circle $\chi$ we have the following possibilities (see \cite{BolsinovF04,Marco09,BolsinovBM10,Marco12}):
\begin{itemize}
\item[(i)] If the index of $\chi$ is $0$ or $2$, then there exists a neighborhood of $\chi$ that is foliated by regular $2$-dimensional tori of the type $I_1=\textrm{const.}$, and  the system of coordinates $(I_1,\phi_1,\phi_2)$ is defined in a cylindrical neighborhood of $\chi$. The intersections of these tori with a surface of section $\mathcal{D}_{0,c}$ are circles surrounding $\chi\cap\mathcal{D}_{0,c}$.
\item[(ii)] If the index of $\chi$ is $1$, then there exists a  neighborhood of $\chi$ that is locally the product between $\chi$ and a `cross' composed of separatrices of the orbit. The trajectory has
an orientable separatrix diagram or a non-orientable separatrix
diagram (see \cite{Fomenko86}). The whole connected component of the critical level of $K$ containing $\chi$
is a finite union of critical circles and cylinders  whose boundary is
either made of one or two critical circles; all these critical circles  have index~$1$. This structure is referred to in \cite{Marco12} as a polycycle. The intersection between a polycycle  and  $\mathcal{D}_{0,c}$  determines separatrices for the dynamics of the Poincar\'e map  $f_{0,c}$.
\end{itemize}

\begin{figure}
\centering
\includegraphics[width=0.5\textwidth]{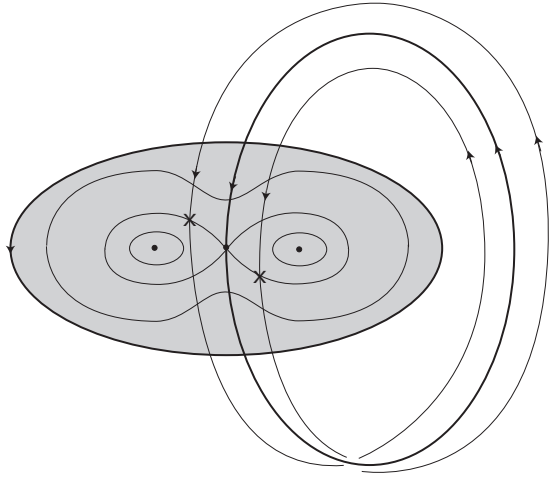}
\caption{Example of a possible topology for the unperturbed system.}
\label{fig:example}
\end{figure}

Consequently, the family of foliations by $2$-dimensional tori of the action-angle domains $\Omega^j_{0,c}$  determine a partition of the surface of section $\mathcal{D}_{0,c}$  into a finite collection of annuli (including punctured disks) $\widehat\Omega_{0,c}^j$, $j=1,\ldots,k$. 
The boundaries of these annuli consist  of points, circles, and separatrices. As it was explained in the proof of Proposition \ref{prop:kam}, each annulus is foliated by invariant circles of the type $J=\textrm{const.}$, and  the Poincar\'e map $f_{0,c}$  satisfies a monotone twist condition on each annulus. An example illustrating a possible topology of the intersection between these objects -- $2$-dimensional tori foliations, critical circles, and polycycles -- with the disk $\mathcal{D}_{0,c}$ is shown in  Fig.~\ref{fig:example}.

\subsubsection{The geometry of the perturbed system}\label{sec:geompert}
We now look into the effect of a small perturbation on the foliations by $2$-dimensional tori, the critical circles, and the polycycles that organize  $\Lambda_{0,c}$.

For $0\leq\eps<\eps_1$ and $c_0\leq c\leq c_1$, the domains $\Omega_{0,c}^j$, $j=1,\ldots,k$,  in $\Lambda_{0,c}$  can be continued to domains $\Omega_{\eps,c}^j$ in $\Lambda_{\eps,c}$, which can still be parametrized by one action and two angle coordinates $(I_1,\phi_1,\phi_2)$.

We make the following additional assumption on the perturbation $H_1$:
\begin{itemize}
\item [(A3)] For every $\eps\in(0,\eps_1)$ and every $c\in[c_0,c_1]$ the following hold true:
\begin{itemize}
\item[(i)] For each component  of the boundary $\partial \Omega_{0,c}^j$ of an action-angle domain  in $\Lambda_{0,c}$ which is of the type $W^u_{\mid \Lambda_{0,c}}(\chi_{0,c})=W^s_{\mid\Lambda_{0,c}}(\chi_{0,c})$, with $\chi_{0,c}$ a critical circle of hyperbolic type, for the perturbed system we have $W^u_{\mid\Lambda_{\eps,c}}(\chi_{\eps,c})$ intersects transversally $W^s_{\mid\Lambda_{\eps,c}}(\chi_{\eps,c})$ in $\Lambda_{\eps,c}$, where $\chi_{\eps,c}$ is the hyperbolic circle corresponding to $\chi_{0,c}$ after perturbation.
    The above stable and unstable manifolds are $2$-dimensional hyperbolic manifolds for the Hamiltonian flow restricted to $\Lambda_{0,c}$, and, respectively $\Lambda_{\eps,c}$.
    \item[(ii)] There exists
a collection of homoclinic channels $\Gamma_{\eps,c}^i$, $i=1,\ldots, m$, and corresponding scattering maps $S_{\eps,c}^i:U^i_{\eps,c}\to V^i_{\eps,c}$ associated to $\Gamma_{\eps,c}^i$, with $U_{\eps,c}^i,V_{\eps,c}^i\subseteq \Lambda_{\eps,c}$, such that, for every $j\in\{1,\ldots,k\}$, each level set of $I_1$ in $\Omega^j_{\eps,c}$ intersects some $U^i_{\eps,c}$.
Moreover, for every $j$, and every action level set $\{I_1=I^0_1\}$ in $\Omega^j_{\eps,c}$, there exist $i$ and  a point $x_d\in U^i_{\eps,c}\cap\{I_1=I_1^0\}$ such that $I_1(S_{\eps,c}^i(x_1))=I_1^d<I_1^0$, and also a point $x_u\in U_{\eps,c}^i\cap\{I_1=I_1^0\}$ such that $I_1(S_{\eps,c}^i(x_2))=I_1^u>I_1^0$.
\item [(iii)] For every critical circle $\chi_{\eps,c}$ that is of hyperbolic type, there exist  scattering map  $S^i_{\eps,c}:U^i_{\eps,c}\to V^i_{\eps,c}$,  with $U^i_{\eps,c}$ intersecting $W^{u,s}_{\mid\Lambda_{\eps,c}}(\chi_{\eps,c})$, that move points from either side of $W^{u,s}_{\mid\Lambda_{\eps,c}}(\chi_{\eps,c})$ to the opposite side in $\Lambda_{\eps,c}$.
\end{itemize}
\end{itemize}

We will show in Subsection \ref{subsub:generic} that  the assumption (A3) is satisfied by an open and dense set of perturbations $H_1$.
Assumption (A3-i) implies that the separatrices of the action-angle domains that correspond to stable and unstable manifolds of hyperbolic circles are destroyed by the perturbation, yielding transverse homoclinic intersections.
Assumption (A3-ii) implies that there exists a collection of scattering maps that is rich enough so that it covers the whole action range inside each  maximal action-angle domain in $\Lambda_{\eps,c}$,  as well as the stable/unstable manifolds  of hyperbolic invariant circles that separate these domains.
Moreover, it says that inside such a domain one can use one of the  scattering maps to  increase/decrease the action coordinate $I_1$.
Assumption (A3-iii) also says that the scattering maps can be used to cross the stable/unstable manifolds of the hyperbolic circles from one side to the other.

Relative to the surface of section $\mathcal{D}_{\eps,c}$, assumption (A3-1) implies that the separatrices that were bounding different action-angle domains $\hat\Omega^j_{0,c}$ are destroyed by the perturbation, and they give rise to hyperbolic periodic orbits  whose stable and unstable manifolds intersect transversally.
Assume that  $O(p)=\{p_1,\ldots, p_m\}$ is a hyperbolic periodic orbit for $f_{\eps,c}$, with $W^u_{\mid \mathcal{D}_{\eps,c}}(O(p))$ transverse to $W^s_{\mid \mathcal{D}_{\eps,c}}(O(p))$  in $\mathcal{D}_{\eps,c}$. These stable and unstable manifolds are $1$-dimensional hyperbolic invariant manifolds for the map $f_{\eps,c}$ on $\mathcal{D}_{\eps,c}$.

Proposition \ref{prop:kam} implies that  inside each domain $\widehat\Omega_{\eps,c}^j$, $j=1,\ldots,k$. 
there exists a family of $1$-dimensional tori $\{\mathcal{T}_\alpha\}_{\alpha\in \mathcal{A}_j}$ that are invariant under $f_{\eps,c}$. Inside each annulus  bounded by a pair of invariant tori $\mathcal{T}_{\alpha_1}$, $\mathcal{T}_{\alpha_2}$, for each rotation number in the corresponding rotation interval, there exists an Aubry-Mather set of that rotation number. A collection of invariant sets $\{\mathcal{R}_1,\ldots, \mathcal{R}_n\}$, where each $\mathcal{R}_i$ is either an invariant $1$-dimensional torus or an Aubry-Mather set, possibly lying in different regions  $\widehat\Omega_{\eps,c}^j$, is given. To show the existence of a trajectory of the Hamiltonian flow $\phi^t_\eps$ that visits these sets, we have to combine the inner dynamics, given by the map $f_{\eps,c}$ on $\mathcal{D}_{\eps,c}$, with the outer dynamics, along homoclinic trajectories corresponding to the homoclinic manifolds $\Gamma^i_{\eps,c}$, $i=1,\ldots,m$.

We need to address several issues. First, the inner dynamics is given in terms of the map $f_{\eps,c}$  on a surface of section on $\mathcal{D}_{\eps,c}$, and the outer dynamics is given in terms of the flow $\phi^t_\eps$ in $\mathbb{R}^6$; we have to translate all information in the language of maps, and reduce the problem to $\mathcal{D}_{\eps,c}$. Second, we have to combine  the inner and outer dynamics, relative to  $\mathcal{D}_{\eps,c}$, to obtain orbits of the discrete dynamical system that visit the prescribed collection of invariant sets. Then, we can conclude that there also exist trajectories of the flow that visit the prescribed collection of invariant sets.

\subsubsection{Reduction to a discrete dynamical system}\label{sec:diff}

We now describe how to combine the outer dynamics, given by homoclinic trajectories and corresponding  scattering maps, the inner dynamics, given by the restriction of the flow to $\Lambda_{\eps,c}$, and the reduced dynamics, given by the return map $f_{\eps,c}$ to the surface of section $\mathcal{D}_{\eps,c}$, in order to construct orbits with the desired characteristics. The main part of the argument relies heavily on the constructions from \cite{GideaR12}; hence we will not repeat some of the technical details, but focus on how to utilize these constructions in order to prove the main result. A description of the ingredients of \cite{GideaR12} is given in Appendix \ref{app:topmet}.

First, we show how to produce trajectories that stay close enough to $\Lambda_{\eps,c}$ for an arbitrarily long time.
We use the linearization of the flow $\phi^t_\eps$ restricted to the energy hypersurface $\{H_\eps=c\}$ in a neighborhood of the normally hyperbolic invariant manifold $\Lambda_{\eps,c}$; see Appendix \ref{section:linearization}.
Let $\mathcal{U}_{\eps,c}$ be a   $(\delta/2)$-neighborhood of $\Lambda_{\eps,c}$ in the energy hypersurface where the flow can be linearized.
Given any time $T>0$, each point in the set $\mathcal{U}_{\eps,c}^{T}=\bigcap_{t\in[0,T]}\phi^{-t}_{\eps}(\mathcal{U}_{\eps,c})$ stays in $\mathcal{U}_{\eps,c}$ for  at least a time $T>0$. The  trajectory can be chosen so that, after it spends at least a time $T\geq 0$ in  $\mathcal{U}_{\eps,c}$, it  follows  the homoclinic orbit corresponding to a specific scattering map $S^i_{\eps,c}$, and then re-enters a neighborhood  $\mathcal{U}_{\eps,c}^{T'}$ corresponding to some other prescribed time $T'>0$.

We discretize the dynamics by considering the time-$1$ map  $F_{\eps}$ of the flow $\phi^t_\eps$, which is defined on the energy hypersurface.
The above considerations for the flow dynamics relative to $\mathcal{U}_{\eps,c}^{T}$ remain valid for the time-$1$ map dynamics (we can choose $T \in\mathbb{Z}^+$).
We also remark that the scattering map for the time-$1$ map coincides with the scattering map for the flow \cite{DelshamsLS08a}.
Hence, we will use the same notation $S^i_{\eps,c}$ to refer to the scattering map for the flow $\phi^t_\eps$ and the scattering map for the time discretization of the flow $F_{\eps}$.

We can also define a scattering map for the return map $f_{\eps,c}$ to $\mathcal{D}_{\eps,c}$.   For each scattering map $S^i_{\eps,c}$ for the flow/time discretization, there exists a scattering map $\widehat S^i_{\eps,c}$ that enjoys similar properties. See \cite{DelshamsGR13}.

We describe a procedure  to combine the dynamics of the time-$1$ map $F_{\eps}$ restricted to $\Lambda_{\eps,c}$  with the dynamics of the return map $f_{\eps,c}$ on $\mathcal{D}_{\eps,c}$. Consider a point $z\in\mathcal{D}_{\eps,c}$ and its image $f^{n}_{\eps,c}(z)$ under some power of $f_{\eps,c}$.
Then $f^{n}_{\eps,c}(z)=\phi^{t(z,n)}_\eps(z)$ for some time $t(z,n)\in\mathbb{R}$. Consider the orbit   of $z$ under the time-discretization map $F_{\eps}$.  To each $n\in\mathbb{Z}$ there exists a unique $n'(z,n)\in\mathbb{Z}$ such that $n'(z,n)\leq t(z,n)<n'(z,n)+1$. The points $f^n_{\eps,c}(z),F^{n'(z,n)}_{\eps}(z)$ lie on the same  trajectory $\phi^t_\eps(z)$. In this way, for each point of an orbit $f^n_{\eps,c}(z)$ of $f^n_{\eps,c}$ we can associate, in a canonical way, a point of an orbit $F^{n'}_{\eps}(z)$ of $F^{n'}_{\eps}$. 

We will use this idea in the following way. We will  combine  the dynamics of $f_{\eps,c}$ and of $\widehat S^i_{\eps,c}$  on $\mathcal{D}_{\eps,c}$ to follow the prescribed invariant objects $\{\mathcal{R}_1,\ldots, \mathcal{R}_n\}$. More precisely, we will construct a sequence of
$2$-dimensional windows  in $\mathcal{D}_{\eps,c}$, with some of the windows in the sequence near the sets $\mathcal{R}_i$, such that any two consecutive windows in the sequence are correctly aligned either under some power of $f_{\eps,c}$,  or under one of the scattering maps $\widehat S^i_{\eps,c}$, $i=1,\ldots,m$. Then, we will use these windows to construct another sequence of $3$-dimensional windows in $\Lambda_{\eps,c}$, that are correctly aligned either under some power of $F_{\eps}$,  or under one of the scattering maps $S^i_{\eps,c}$, $i=1,\ldots,m$. Then we will invoke Lemma \ref{lem:shadowing2} to show that there exists orbits of $F_{\eps}$ in $\mathbb{R}^6$ that follow this latter sequence of windows. We will then conclude that there exist orbits of the flow $\phi^t_\eps$ that follow these windows.

\subsubsection{Construction of a sequence of correctly aligned $2$-dimensional windows}\label{section:2win}
Given the collection  of invariant sets $\{\mathcal{R}_1,\ldots,\mathcal{R}_n\}$ in $\mathcal{D}_{\eps,c}$, we would like to construct orbits that visit, in the consecutive order, all the sets $\mathcal{R}_i$ that lie within the same  domain    $\widehat\Omega^j_{\eps,c}$, then move to an adjacent domain and visit all  the sets $\mathcal{R}_{i'}$ that lie within  that domain, and so on, until all sets   $\mathcal{R}_i$ have been visited.

Suppose that a sub-collection $\{\mathcal{R}_{j_1},\ldots,\mathcal{R}_{j_n}\}$ of  $\{\mathcal{R}_1,\ldots,\mathcal{R}_n\}$  is contained in some $\widehat\Omega^j_{\eps,c}$. 
Recall that the  sets   $\mathcal{R}_i$ are either $1$-dimensional KAM tori or Aubry-Mather sets.
We consider the entire family of KAM tori $\{\mathcal{T}_{\alpha}\}_{\alpha\in\mathcal{A}_j}$ provided by Proposition \ref{prop:kam}, which includes among its members the sets $\mathcal{R}_i$ from the above sub-collection that are KAM tori. Recall that the family of scattering maps $S^i_{\eps,c}$, $i=1,\ldots,m$, yields a family of scattering maps $\widehat S^i_{\eps,c}$,  with $\widehat S^i_{\eps,c}: \widehat U^i_{\eps,c}\to \widehat V^i_{\eps,c}$, where $\widehat U^i_{\eps,c}, \widehat V^i_{\eps,c}$ are open subsets of  $\mathcal{D}_{\eps,c}$, for $i=1,\ldots,m$.

Assumption (A3-ii) implies that each level set $J=J^0$ in $\widehat\Omega^j_{\eps,c}$ intersects the domain of some scattering map $\widehat S^i_{\eps,c}$, and, moreover, there exist   points $\hat x^d$ on that level set such that $J(\widehat S^i_{\eps,c}(\hat x^d))=J^d<J^0$, and there exist points $\hat x^u$ on that level set such that $J(\widehat S^i_{\eps,c}(\hat x^u))=J^u>J^0$.

Using these scattering maps and the inner dynamics  we can construct in $\widehat\Omega^j_{\eps,c}$  a new sequence consisting of $1$-dimensional KAM tori and essential circles that are not necessarily KAM tori,  which we denote by $\{\mathcal{T}_{\tau}\}_{\tau=1,\ldots, \tau_j}$, for some integer $\tau_j>0$, and is characterized by the following properties:
\begin{itemize}
\item The sequence $\{\mathcal{T}_{\tau}\}_{\tau=1,\ldots,\tau_j}$ consists of subsequences  $\{\mathcal{T}_{\tau'}, \mathcal{T}_{\tau'+1},\ldots, \mathcal{T}_{\tau'+l}\}\subseteq\{\mathcal{T}_{\tau}\}_{\tau=1,
\ldots, \tau_j}$  that form {\em transition chains of tori}, alternating with {\em Birkhoff Zones of Instability}, between the tori $\mathcal{T}_{\tau'+l}$ and $\mathcal{T}_{\tau'+l+1}$:
     \begin{itemize}
\item the transition chain  is characterized by the fact that for  each $t=0,\ldots, l-1$ there is an $i\in\{1,\ldots,m\}$ such that  $\widehat S^i_{\eps,c}(\mathcal{T}_{\tau'+t})$  topologically crosses  $ \mathcal{T}_{\tau'+t+1}$  (for the definition of topological crossing see,  e.g.,  \cite{GideaR03});
\item the Birkhoff Zones of Instability between $\mathcal{T}_{\tau'+l}$ and $\mathcal{T}_{\tau'+l+1}$ is characterized by the fact that the region in
$\widehat\Omega^j_\eps$ between $\mathcal{T}_{\tau'+l}$ and $\mathcal{T}_{\tau'+l+1}$ does not contain any essential invariant circle in its
interior;
\end{itemize}
\item Each torus $\mathcal{T}_\tau$ that is  not an  end torus of a transition chain is $C^1$-smooth, and each torus $\mathcal{T}_\tau$ that is  an  end torus of a transition chain is only Lipschitz continuous;
\item Each torus $\mathcal{T}_\tau$ is  the limit  of some other invariant tori (not necessarily from this sequence), i.e., there exists a sequence of $1$-dimensional invariant tori $(\mathscr{T}_n)_{n\in\mathbb{N}}$ that approaches $\mathcal{T}_\tau$ in the $C^0$-topology:
\begin{itemize}
\item if $\mathcal{T}_\tau$  is an end torus of a transition chain, then it can be approximated by  some other invariant tori from only one side;
\item if $\mathcal{T}_\tau$  is not an  end torus of a transition chain, then it can be approximated by  some other invariant tori from both sides;
\end{itemize}
\item Every   $\mathcal{R}_i$ from the sub-collection $\{\mathcal{R}_{j_1},\ldots,\mathcal{R}_{j_n}\}$
that is an invariant torus, is a member of one of the transition chains from above, and every  $\mathcal{R}_i$ that is an Aubry-Mather set lies in one of the Birkhoff Zones of Instability from above.
\end{itemize}

We now invoke  the construction of correctly aligned windows from \cite{GideaR12}.
That construction yields  a sequence  of $2$-dimensional windows in $\widehat\Omega^j_{\eps,c}$, which we denote by
\[\widehat R_0,\widehat R_1,\ldots,\widehat R_{p},\]
such that:
\begin{itemize}
\item For every  $q=0,\ldots,p-1$, $\widehat R_q$ is correctly aligned with $\widehat R_{q+1}$ under some power $f^{n_q}_{\eps,c}$ of $f_{\eps,c}$, or under some scattering map $\widehat S^i_{\eps,c}$,  where $i=1,\ldots,m$;
\item For every set $\mathcal{R}_i$ there exists a window $\widehat R_q$ in the above sequence  that lies within a $(\delta/2)$-neighborhood of $\mathcal{R}_i$.
\end{itemize}

By the assumptions (A3-i) and (A3-iii), there are orbits that  go from one   domain $\widehat\Omega_{\eps,c}^j$ to another. We now construct correctly aligned windows along such orbits are link them with the windows along the sequences of tori constructed above. We do the following. Assume that $\chi_{\eps,c}$ is a critical circle of hyperbolic type with $W^u_{\mid \Lambda_{\eps,c}}(\chi_{\eps,c})$ intersecting transversally $W^s_{\mid \Lambda_{\eps,c}}(\chi_{\eps,c})$, as in (A3-i).  Let the intersection of $\chi_{\eps,c}$ with $\mathcal{D}_{\eps,c}$ be the hyperbolic periodic orbit $O(p)$ for $f_{\eps,c}$, and let  $W^u_{\mathcal{D}_{\eps,c}}(O(p))$ and $W^s_{\mathcal{D}_{\eps,c}}(O(p))$ be the corresponding hyperbolic invariant manifolds in $\mathcal{D}_{\eps,c}$; they do also intersect transversally in $\mathcal{D}_{\eps,c}$. Let $\mathcal{T}$ be the last invariant circle in $\widehat\Omega^j_{\eps,c}$, defined by the property that there is no essential  invariant circle between $\mathcal{T}$ and $W^u_{\mathcal{D}_{\eps,c}}(O(p))\cup W^s_{\mathcal{D}_{\eps,c}}(O(p))$. We  extend the sequence of tori $\{\mathcal{T}_{\tau}\}_{\tau=1,\ldots,\tau_j}$ to include the torus $\mathcal{T}$, and  extend the construction of correctly aligned windows to attain $\mathcal{T}$.
The region between  $\mathcal{T}$ and $W^u_{\mathcal{D}_{\eps,c}}(O(p))\cup W^s_{\mathcal{D}_{\eps,c}}(O(p))$ is also a BZI.  We can again call the arguments in \cite{GideaR12} to construct a window $\widehat R$ by $\mathcal{T}$ which is correctly aligned under some power of $f_{\eps,c}$ with a window $\widehat R^s$ by $W^s_{\mathcal{D}_{\eps,c}}(O(p)))$. By the standard construction of the correctly aligned windows by a hyperbolic periodic orbit, we can construct another window $\widehat R^u$ by $W^u_{\mathcal{D}_{\eps,c}}(O(p))$, such that $\widehat R^s$ is correctly aligned with $\widehat R^u$ under some power of $f_{\eps,c}$. By (A3-iii) we can also  move from one side  of $W^{u,s}_{\mathcal{D}_{\eps,c}}(O(p))$ to another.
In this way we  obtain  that the construction of correctly aligned can also be extended across different domains $\widehat\Omega^j_{\eps,c}$.

At the next step, we use the above $2$-dimensional windows in $\mathcal{D}_{\eps,c}$, which are correctly aligned under $f_{\eps,c}$, to construct $3$-dimensional windows in $\Lambda_{\eps,c}$, that are correctly aligned under $F_{\eps}$.

\subsubsection{Construction of a sequence of correctly aligned $3$-dimensional windows}\label{section:3win}

For each $2$-dimensional window $\widehat R_q\subseteq \mathcal{D}_{\eps,c}$ we construct a  $3$-dimensional windows $R_q\subseteq\Lambda_{\eps,c}$,  with the property that $R_q$ is correctly aligned with $R_{q+1}$ either under some suitable power of $F_{\eps}$ or under some scattering map $S^j_{\eps,c}$, $j=1,\ldots,k$. 
Note that  the flow lines of $\phi^t_\eps$ restricted to $\Lambda_{\eps,c}$ which are passing through   $\widehat R_q$  intersect   $\widehat R_{q+1}$ along $f^{n_q}_{\eps,c}(\widehat R_q)\cap \widehat R_{q+1}$.

We start with the construction of $R_0$ and continue the construction inductively. Suppose  that $\widehat R_0$ is correctly aligned with $\widehat R_1$ under $f^{n_0}_{\eps,c}$.  Choose a  point $z\in \widehat R_0$  and let $n'(n_0):=n'_0$ be the unique integer with $n'_0\leq t(z,n_0)<n'_0+1$. Since $F_{\eps}$ is a diffeomorphism, $F^{n'_0}_{\eps}(\mathcal{D}_{\eps,c})$ is a diffeomorphic copy of $\mathcal{D}_{\eps,c}$ and  is itself a global surface of section (sharing the same boundary circle with $\mathcal{D}_{\eps,c}$).
Define $\widecheck R_1$ as the the set of points of first intersection between the flow lines through $\widehat R_1$ and  $F^{n'_0}_{\eps}(\mathcal{D}_{\eps,c})$.
We define the window $R_0$  of the type
\[R_0=\bigcup_{z\in\widehat{R}_0}\phi^{[t_1(0,z),t_2(0,z)]}_\eps(z),\] for some $t_1(0,z),t_2(0,z)\in\mathbb{R}$, with the property that $\widehat R_0= \bigcup_{z\in\widehat{R}_0}\phi^{t_0(0,z)}_\eps(z)$ for some $t_0(0,z)\in (t_1(0,z),t_2(0,z))$.
We let
\[R_0^{\rm ex} =\bigcup_{z\in\widehat{R}_0^{\rm ex}}\phi^{[t_1(0,z),t_2(0,z)]}_\eps(z) \cup \bigcup_{z\in\widehat{R}_0}\phi^{\partial[t_1(0,z),t_2(0,z)]}_\eps(z),\]
be the exit set of $R_0$.

Then we define
\[R _{1}=\bigcup_{z\in\widecheck{R}_{1}}\phi^{[t_1( 1,z),t_2( 1,z)]}_\eps(z),\] for some $t_1( 1,z),t_2( 1,z)\in\mathbb{R}$, with the property that $\widecheck R_{ 1}= \bigcup_{z\in\widecheck{R}_{ 1}}\phi^{ t_0( 1,z) }_\eps(z)$ for some $t_0( 1,z)\in (t_1(1,z),t_2(1,z))$.
We let
\[R_1^{\rm ex} =\bigcup_{z\in\widecheck{R}_1^{\rm ex}}\phi^{[t_1(1,z),t_2(1,z)]}_\eps(z) \cup \bigcup_{z\in\widecheck{R}_1}\phi^{\partial[t_1(1,z),t_2(1,z)]}_\eps(z).\]

\begin{figure}
\centering
\includegraphics[width=0.8\textwidth]{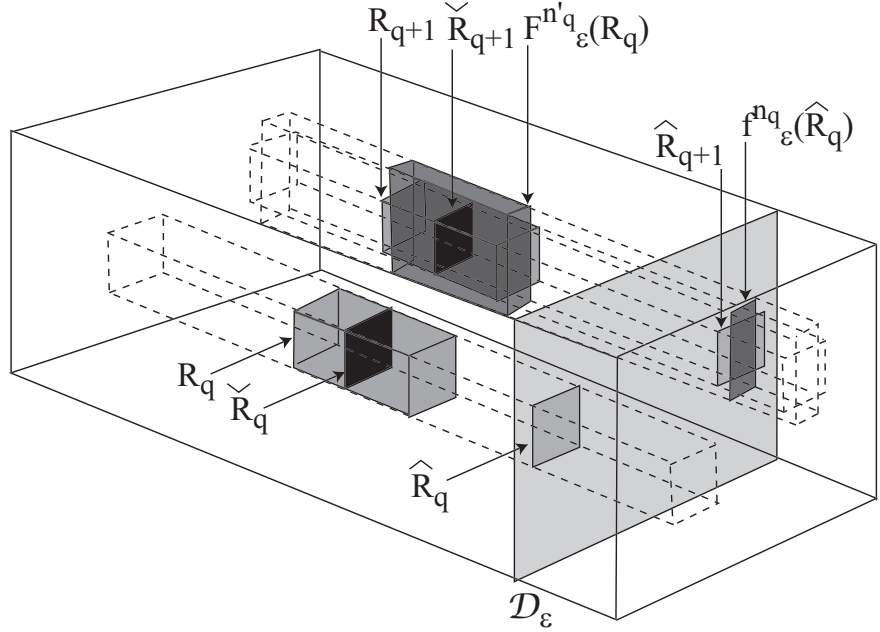}
\caption{Windows correctly aligned by the return map and windows correctly aligned by the time-$1$ map.}
\label{fig:returnmapdynamics}
\end{figure}

It is immediate that for some suitable choices of the functions $t_1(0,z) ,t_2(0,z) ,t_0(0,z)$, $t_1(1,z) ,t_2(1,z) ,t_0(1,z)$,  we get that
$R_0$ is correctly aligned with $R_{1}$ under $F^{n'_0}_{\eps,c}$. More precisely, we require that $t_1(0,z)+n'_0< t_1(1,z)<t_2(1,z)<t_2(0,z)+n'_0$.
See Fig. \ref{fig:returnmapdynamics}.

In a similar manner,  given a pair of windows $\widehat {R}_q, \widehat {R}_{q+1}$ in $\mathcal{D}_{\eps,c}$, with  $\widehat {R}_q$ correctly aligned with $\widehat {R}_{q+1}$  under $\widehat S^j_{\eps,c}$,
we can construct a  pair of windows $ {R}_q,  {R}_{q+1}$ in $\Lambda_{\eps,c}$,   with  ${R}_q$ correctly aligned with  ${R}_{q+1}$ under $S^j_{\eps,c}$.

The construction is continued inductively, providing a sequence
of windows
\[R_0, R_1,\ldots, R_{p},\] such that $R_q$ is correctly aligned with $R_{q+1}$, either under some power $F^{n'_q}_{\eps,c}$ of $F_{\eps}$, or under some scattering map $S^j_{\eps,c}$, $j=1,\ldots,k$.
 It is implicit in this construction
that the trajectories of the flow $\phi^t_\eps$ that that pass through the windows $R_q$, $q=0,\dots,p$,  also pass $(\delta/2)$-close to each of the prescribed sets $\mathcal{R}_i$, $i=1,\ldots,m$.

We apply Lemma \ref{lem:shadowing2} from Appendix \ref{app:windows}, implying  that for correctly aligned windows as the ones constructed above there exist flow trajectories that visit $(\delta/2)$-neighborhood of these windows, in the prescribed order. Hence these flow trajectories pass $\delta$-close to each of the prescribed sets $\mathcal{R}_i$, $i=1,\ldots,m$.
This completes the argument.

\begin{rem}
We remark that we can also `thicken' the windows $\widehat R_q$ to $3$-dimensional windows $\widehat{\widehat {R}}_q$ of the type
\begin{equation*}
\begin{split} \widehat{\widehat {R}}_q=&\bigcup_{z\in\widehat{R}_q}\phi^{[t_1(q,z),t_2(q,z)]}_{\eps}(z), \\
 \widehat{\widehat {R}}_q=&\bigcup_{z\in\widehat{R}_q^{\rm ex}}\phi^{[t_1(q,z),t_2(q,z)]}_\eps(z) \\
 &\cup \bigcup_{z\in\widehat{R}_q}\phi^{\partial[t_1(q,z),t_2(q,z)]}_\eps(z),
\end{split}
\end{equation*}
for some suitable $t_1(q,z),t_2(q,z)\in\mathbb{R}$, with the property that $\widehat R_q= \bigcup_{z\in\widehat{R}_q}\phi^{ t_0(q,z) }_\eps(z)$ for some $t_0(q,z)\in (t_1(q,z),t_2(q,z))$.
Then the relation between the windows $R_q$ and the windows $\widehat{\widehat {R}}_q$ is that each window $R_q$ is correctly aligned with the corresponding $\widehat{\widehat {R}}_q$ under some suitable translation mapping along the trajectories of the flow $\phi^t_\eps$ through $\widehat R_q$.
\end{rem}

\begin{rem}
We  have used the topological method of correctly aligned windows to show the existence of orbits that shadow Aubry-Mather sets
lying in different maximal action-angle domains.  The existence of orbits that shadow  Aubry-Mather sets within a single action-angle domain has been established through variational methods in \cite{Mather91}. This approach has been applied to the diffusion problem in \cite{ChengY2004,ChengY2009}. The orbits found in these papers are minimal. As concatenations of minimal orbit does not necessarily yield minimal orbits, the variational approach does not seem to be suitable to find orbits traveling across different action-angle domains.
\end{rem}

\begin{rem}
One of the main ingredients to achieve global diffusion in this problem is the combination of  multiple dynamics, i.e., of  the inner dynamics given by the restriction of the Hamiltonian flow to $\Lambda_{\eps,c}$ and of the outer dynamics given by several scattering maps. This is closely related to the idea of polysystems, also known as iterated  function systems, which has enjoyed  a recent development on its own (see, for example  \cite{Moeckel99,Hutchinson81,MihailM08,MiculescuM08,MiculescuM13}).
\end{rem}

\subsubsection{Genericity of the perturbation}\label{subsub:generic}
In this section we sketch out the argument that there exists an `generic' set of Hamiltonians $H_1$ for which condition (A3) from Subsection \ref{sec:geompert} is satisfied. We follow the ideas in \cite{DelshamsLS2006}, where detailed arguments can be found. We can restrict the domain of the Hamiltonians $H_\eps$, $\eps\in[0,\eps_1]$ to some conveniently large compact set $\mathscr{K}\subseteq \mathbb{R}^6$. By a generic set of Hamiltonians $H_1$ we mean an open and dense set of Hamiltonians in $C^r(\mathscr{K})$.

We assume that a Hamiltonian $H_1\in C^r(\mathscr{K})$ is given, and we show that by an arbitrarily small perturbation of  $H_1$ in $C^r(\mathscr{K})$ the resulting Hamiltonian $H_\eps=H_0+\eps H_1$ satisfies the condition (A3).

By assumption (A2) the Hamiltonian $H_{01}$ has two homoclinic orbits $\gamma^\pm$   to $(0,0)$. We denote by $ (p^\pm_3(\tau),q^\pm_3 (\tau))$, $\tau\in\mathbb{R}$, some parametrizations of these homoclinic orbits with
$(p^\pm_3(0),q^\pm_3 (0))=(0,0)$, where $(0,0)$ is the saddle equilibrium point of the pendulum-like system. These homoclinic
orbits determine two corresponding branches of the stable and unstable manifolds of
$\Lambda_{0,c}$, $W^{+,u}(\Lambda_{0,c})=W^{+,s}(\Lambda_{0,c})=\bigcup_{z\in\Lambda_{0,c}}\gamma^+$, and
and $W^{-,u}(\Lambda_0)=W^{-,s}(\Lambda_0)=\bigcup_{z\in\Lambda_{0,c}}\gamma^-$.

We want to measure the splitting of the stable and unstable manifolds of $\Lambda_\eps$ when $\eps\neq 0$ is sufficiently small. For this, we define the  Melnikov
potential for each of the homoclinic orbits $\gamma^\pm$  by
\begin{eqnarray}\label{eqn:melnikov}
\mathcal{M}^\pm(z,\tau)=&\displaystyle \int
_{-\infty}^{\infty}&\left
[H_1(p_1(t),q_1(t),p_2(t),q_2(t), p^\pm_3(\tau+t),q^\pm_3(\tau+t))
\right .\\
& &\left .-H_1(p_1(t),q_1(t),p_2(t),q_2(t), 0,0)\right ]dt,
\end{eqnarray}
where $z=(p_1,q_1,p_2,q_2)\in \Lambda_{0,c}$, and $(p_1(t),q_1(t),p_2(t),q_2(t))$ denotes the trajectory   $\phi^t_0(z)$ on $\Lambda_0$ for the Hamiltonian  $H_{00}$.

If we restrict to a domain $\Omega^j_{0,c}\subseteq \Lambda_{0,c}$, with corresponding coordinates $(I_1,\phi_1,\phi_2)$, with respect to these coordinates we have
\begin{eqnarray*}
\mathcal{M}^\pm(I_1,\phi_1,\phi_2,\tau)=&\displaystyle \int
_{-\infty}^{\infty}&\left
[H_1(I_1,\phi_1+tI_1, I_2,\phi_2+tI_1,  p^\pm_3(\tau+t),q^\pm_3(\tau+t))\right .\\
& &\left .-H_1(I_1,\phi_1+tI_1, I_2,\phi_2+tI_1,0,0)\right ]d t,
\end{eqnarray*}
where the value $I_2$ is implicitly determined by $(I_1,\phi_1,\phi_2)$.
To measure the splitting of the stable and unstable manifolds of $\Lambda_{\eps,c}$, we take a local section $\Sigma$ to $\gamma^\pm$ through $\gamma^\pm(\tau)$ and we measure the distance between the intersection points of $W^{\pm,u}(k_\eps(z))$ and
 $W^{\pm,s}(k_\eps(z)) $ with $\Sigma$, which turns out to be given by
\begin{eqnarray*}
-\eps\frac{d}{d\tau}\mathcal{M}^\pm(I_1,\phi_1,\phi_2,\tau)+O(\eps^{1+\varrho}),
\end{eqnarray*}
for some $\varrho>0$.
The existence of a transverse intersection  of the stable and unstable manifolds $W^{\pm,s}(\Lambda_{\eps,c})$, $W^{\pm,u}(\Lambda_{\eps,c})$ is guaranteed provided that the function  $\tau\mapsto \mathcal{M}^\pm(I_1,\phi_1,\phi_2,\tau)$ has a non-degenerate critical point $\tau^*(I_1,\phi_1,\phi_2)$. Then $(I_1, \phi_1,\phi_2)\mapsto (I_1, \phi_1,\phi_2,\tau^*(I_1,\phi_1,\phi_2))$, with $(I_1,\phi_1,\phi_2)$ in some domain $U_{0,c}^0$ in $\Lambda_{0,c}$ where $\tau^*$ is a non-degenerate critical point of $\tau\mapsto \mathcal{M}^\pm$,  gives   a parametrization of a homoclinic manifold $\Gamma^0_{\eps,c}\subseteq W^{\pm,s}(\Lambda_{\eps,c})\cap W^{\pm,u}(\Lambda_{\eps,c})$. If $\Gamma^0_{\eps,c}$ is chosen small enough so that it is a homoclinic channel, then there is an associated scattering map $S^0_{\eps,c}: U^0_{\eps,c}\to V^0_{\eps,c}$. For each $x\in\Gamma^0_{\eps,c}$, $S^0_{\eps,c}$ assigns to the point $\Omega^u_{\eps,c}(x)\in U^0_{\eps,c}$ the point $\Omega^s_{\eps,c}(x)\in V^0_{\eps,c}$, where $U^0_{\eps,c},V^0_{\eps,c}\subseteq \Lambda_{\eps,c}$. See Appendix \ref{app:scattering}.

The change in the action $I_1$ by the corresponding scattering map $S^0_\eps$  is given by
\[I_1(\Omega^s_{\eps,c}(x))-I_1(\Omega^u_{\eps,c}(x))=-\eps\frac{d}{d\tau}\{I_1,\mathcal{M}^\pm\}(I_1,\phi_1,\phi_2,\tau^*(I_1,\phi_1,\phi_2)),\]
for $x\in\Gamma_\eps^0$.

Fixing a value $I_1=I_1^0$, we can always ensure, by  an arbitrarily small perturbation of   $H_1$, if necessary, that there exists an open neighborhood $U^0_\eps$ of a point $(I^0_1,\phi^0_1,\phi^0_2)\in\Lambda_\eps$ for  which the map $\tau \mapsto \mathcal{M}^\pm(I_1,\phi_1,\phi_2,\tau)$ has a non-degenerate critical point $\tau^*(I_1,\phi_1,\phi_2)$,
for each  $(I_1,\phi_1,\phi_2)\in U^0_{\eps,c}$.  We can choose this perturbation to have compact support in some small tubular neighborhood of a point of $(I_1,\phi_1,\phi_2,\gamma^+(\tau))$ or of $(I_1,\phi_1,\phi_2,\gamma^-(\tau))$. The   homoclinic intersection given by $(I_1,\phi_1,\phi_2,\tau^*(I_1,\phi_1,\phi_2))$, with $(I_1,\phi_1,\phi_2)\in U^0_{\eps,c}$, is denoted by $\Gamma^0_{\eps,c}$.

Since in the unperturbed system we have available two homoclinic orbits $\gamma^\pm$ that are geometrically different, we can produce arbitrarily small perturbations of $H_1$ that have mutually disjoint supports,  and obtain scattering maps $S^i_{\eps,c}:U^i_{\eps,c}\to V^i_{\eps,c}$,  with corresponding homoclinic manifolds $\Gamma^i_{\eps,c}$,  such that the union of the domains $U^i_\eps$ covers the whole range of $I_1$. This ensures the genericity of condition (A3-ii) within each angle-action domain $\Omega^j_{\eps,c}$, $j=1,\ldots,k$.

Now we discuss the splitting of the hyperbolic invariant manifolds of the critical circles $\chi_0$. For the Hamiltonian flow of $H_0$, each such circle $\chi_0$ has two dimensional stable and unstable manifolds $W^{\pm,s}(\chi_0)$, $W^{\pm,u}(\chi_0)$, with one hyperbolic direction corresponding to the dynamics of $\phi^t_0$ on $\Lambda_0$, and
another hyperbolic direction corresponding to the separatrix  of the flow of $H_{01}$.  We have that $W^{\pm,s}(\chi_0)= W^{\pm,u}(\chi_0)$. By  an arbitrarily small perturbation of a given  $H_1$, if necessary, we can make that $W^{\pm,s}(\chi_0)$ intersects transversally $W^{\pm,u}(\chi_0)$. This ensures the genericity of the  condition (A3-i), and also  of the condition (A3-iii).

\section*{Acknowledgement} Part of this work has been done while M.G. was a member of the IAS, whose support is kindly acknowledged.  The author would also like to thank to Helmut Hofer, Umberto Hryniewicz, Richard Moeckel, Rafael de la Llave,  and Pedro Salom\~ao for useful suggestions, ideas, and discussions.

\appendix
\section{Background on symplectic dynamics.}\label{section:symplectic}

\subsection{Contact geometry}
Given a  compact, connected, oriented, $3$-dimensional  manifold $M$, a contact form $\lambda$ on $M$ is a $1$-form on $M$ such that $\lambda\wedge d\lambda$ is a volume form on $M$. The contact structure associated to $\lambda$ is the plane bundle  in $TM$ given by $\xi=\textrm{ker}(\lambda)=\{(x,h)\in TM\,|\,\lambda(x)(h)=0\}$.  The restriction $d\lambda_{\mid \xi\oplus\xi}$ defines a symplectic structure on each fiber of $\xi\to M$. The characteristic distribution of $M$ is the $1$-dimensional distribution
\[\mathcal{L}=\{(x,h)\in TM\,|\,x\in M,\, \omega(h,k)=0\textrm { for all } k\in T_xM\}.\]
The corresponding $1$-dimensional foliation is called the characteristic foliation.

The Reeb vector field $X$ associate to $\lambda$ is the vector field on $M$ uniquely defined by $i_{X}(d\lambda)=0$ and  $i_{X}(\lambda)=1$.
The Reeb vector field $X$ spans the characteristic
distribution $\mathcal{L}$ which has the canonical section $X$. The flow lines of $X$ are contained in the leaves of the characteristic foliation.  We have that $TM$ naturally splits as
\[TM=\mathcal{L}\oplus \xi=\mathbb{R}X\oplus \xi.\]

A contact structure $\lambda$ is said to be tight provided that there are no overtwisted disks in $M$, that is, there is no embedded disk $D\subseteq M$ such that $T\partial D\subseteq \xi$ and $T_pD\neq \xi_p$ for all $p\in\partial D$. As an example, the $1$-form on $\mathbb{R}^4$ of coordinates  $(q_1,p_1,q_2,p_2)$,
\[\lambda_0=\frac{1}{2}(q_1dp_1-p_1dq_1+q_2dp_2-p_2dq_2),\] gives a tight contact form on the $3$-dimensional sphere $S^3\subseteq \mathbb{R}^4$ when restricted to $S^3$. By a theorem of Eliashberg \cite{Eliashberg92}, every tight contact form on $S^3$ is diffeomorphic to $g\lambda_0$, for some $C^1$-differentiable $g:S^3\to\mathbb{R}\setminus\{0\}$. The sphere $S^3$ equipped with this distinguished contact
structure is called the tight three-sphere.


If we denote by $\phi^t$ the flow of $X$, we have that $(\phi^t)^*\lambda=\lambda$, and $(D\phi^t)_m:\xi_m\to \xi_{\phi^t(m)}$ is symplectic with respect to $d\lambda$.

\subsection{The Conley-Zehnder index} To each contractible $T$-periodic solution $x(t)$ of the Reeb vector field,  there is assigned the so called Conley-Zehnder index. The Conley-Zehnder index generalizes the usual Morse index for closed geodesics on a
Riemannian manifold. Roughly speaking, the index measures how much neighboring
trajectories of the same energy wind around the orbit.

We now recall  the definition of the index. Assume that $x$ is a $T$-periodic solution, which is contractible. The derivative map $D\phi^t:T_{x(0)}M\to T_{x(t)}M$ maps the contact plane $\xi_{x(0)}$ to $\xi_{x(t)}$ and is symplectic with respect to $d\lambda$.
We assume that $x(t)$ is non-degenerate, meaning that  $D\phi^T$ does not contain $1$ in the spectrum.
Choose a smooth disk $u:D\to M$ s.t. $u(e^{2i\pi t/T})=x(t)$, where $D=\{z\in\mathbb{C}\,|\,|z|\leq 1\}$. Then choose a symplectic trivialization $\beta: u^*\xi\to D\times \mathbb{R}^2$. We associate to $x(t)$  an arc of symplectic matrices $\Phi:[0,T]\to Sp(1)$, where $Sp(1)=\{A\in GL(2)\,|\, A^TJA=J\}$, by
\[ \Phi(t)= \beta(e^{2i\pi t/T})\circ (D\phi^t)_{\mid \xi_{x(t)}}\circ \beta(1)^{-1}.\]
The arc starts at the identity $\Phi(0) = \textrm{Id}$  and ends at $\Phi(T)$, with $\det(\Phi(T)-I)\neq 0$, due to the non-degeneracy condition.  
Take $z\in\mathbb{C}$ and compute the winding number of $\Phi(t)z$,
\[\Delta(z)=\theta(T)-\theta(0)\in\mathbb{R},\]  where $\theta(t)$ is a continuous argument of $\Phi(t)z$, i.e., $\Phi(t)z=r(t)e^{2\pi i \theta(t)}$.
Then define the winding interval of the arc $\Phi(t)$ by
\[I(\Phi)=\{\Delta(z)\,|\,z\in \mathbb{C}\setminus \{0\}\}.\]
Equivalently, we can put $\Phi(t)e^{2\pi i s}=r(t,s)e^{2\pi i\theta(t,s)}$ for all $s\in[0,1]$, where $\theta(0,s)=s$, and define $\Delta(s)=\theta(T,s)-\theta(0,s)=\theta(T,s)-s$, and $I(\Phi)=\{\Delta(s)\,|\,s\in [0,1]\}$.
The length of the winding interval is less than $1/2$. Then the winding interval either lies between two consecutive integers or contains
precisely one integer.

We  define
\[\mu(\Phi)=\left\{
              \begin{array}{ll}
                2k, & \hbox{if $k\in I(\Phi)$;} \\
                2k+1, & \hbox{if $I(\phi)\subset (k,k+1)$.}
              \end{array}
            \right.
 \]

Then we define the Conley-Zehnder index $ \mu(x,T, [u])$ of $(x,T)$ by $ \mu(x,T, [u]) =\mu(\Phi)$. It depends on $x$ and $T$ and on the homotopy class of the  choice of the disk map $u:D\to M$ satisfying $u(e^{2\pi it/T})=x(t)$. In the case when $\pi_2(M)=0$ (e.g., if $M=S^3$),  the index is independent of the choice of the disk map $u$.


\subsection{Existence of global surfaces of section}\label{subsection:global}
Consider $\mathbb{R}^4=\{x=(q_1,p_1,q_2,p_2)
\,|\break \, q_1,p_1,q_2,p_2\in\mathbb{R}\}$ endowed with the standard symplectic form $\omega=\sum_{i=1}^{2}dq_j\wedge dp_j$, and $H:\mathbb{R}^4\to\mathbb{R}$ a $C^r$-differentiable Hamiltonian function. If $c\in\mathbb{R}$ is a regular value for $H$, then $S_c=\{x\,|\, H(x)=c\}$ is a $3$-dimensional manifold invariant under the Hamiltonian flow of $H$. Assume that $S_c$ is compact and connected.

The manifold $S_c$ is said to bound a  strictly convex domain provided that there exists $\delta>0$ such that $D^2H(x)- \delta \cdot \textrm{id}$ is positive definite for all $x\in\mathbb{R}^4$. This is equivalent with the conditions that $W=\{x\,|\, H(x)\leq c\}$  is bounded, and $D^2H(x)(h,h)>0$ for each $x\in W$ and each non-zero vector $h$.

An energy manifold $S_c$ of the Hamiltonian $H$  is said to be of contact type if there exists a one-form $\lambda$ on $S_c$ such that
$d\lambda=-j^*\omega$ and $i_{X_H}(\lambda)\neq 0$ hold on $S_c$, where $j:S_c\to \mathbb{R}^4$ is the inclusion map.


Assume now that $S_c$ is diffeomorphic to $S^3$, that it is of contact type, and that the contact structure is tight. The manifold $S_c$ is said to be  dynamically convex if for every periodic solution $(x,T)$ of the Reeb vector field, we have $\tilde\mu(x,T)\geq 3$.

If $S_c$ is equipped with the contact form ${\lambda_0}_{\mid S_c}$, encloses $0_{\mathbb{R}^4}$ and is strictly convex, then it is dynamically convex. The converse is not true.

The following result provides sufficient conditions for the existence of a disk-like surface of section.   Given a closed $3$-dimensional manifold and a flow $\phi$ with no rest points,  we say that a topologically embedded $2$-dimensional disk $\mathcal{D}$ is a disk-like global surface of section provided that: (i) the boundary $\text{bd}(\mathcal{D})$ is a periodic orbit (called spanning orbit), (ii) the interior of the disk $\text{int}(\mathcal{D})$ is a smooth manifold transverse to the flow, and (iii) every orbit, other than the spanning orbit, intersects $\text{int}(\mathcal{D})$ in forward and backward time.

\begin{thm}[\cite{HZW98}]\label{thm:HWZ98}
Assume that $S_c$ is diffeomorphic to $S^3$, is equipped with a tight contact structure, and is dynamically convex. Then there exits a global disk-like surface of section $\mathcal{D}$ and an associated global return map $f:\text{int}(\mathcal{D})\to \text{int}(\mathcal{D})$ that is smoothly conjugated to an area preserving mapping of the open unit  disk in $\mathbb{R}^2$. The spanning orbit $\chi$ of prime period $T$ has   Conley-Zehnder index $\tilde\mu(\chi,T)=3$.
\end{thm}

We note that a generalization of this result to non-dynamically convex tight contact forms on the three-sphere appears in \cite{HryniewiczS2011}.


\section{Background on the scattering map.}\label{app:scattering} Consider a flow $\Phi: M\times \mathbb{R}\to M$ defined on a manifold $M$ that possesses a normally hyperbolic invariant manifold $\Lambda\subseteq M$.

As the stable and unstable manifolds of $\Lambda$ are foliated by stable and unstable manifolds of points, respectively, for each $x\in W^u(\Lambda)$ there exists a unique $x^u\in\Lambda$ such that
$x\in W^u(x^u)$, and for each   $x\in W^s(\Lambda)$ there exists
a unique $x^s\in\Lambda$ such that $x\in W^s(x^s)$.
We define the  wave maps  $\Omega^s:W^s(\Lambda)\to \Lambda$ by
$\Omega^s(x)=x^u$, and $\Omega^u:W^u(\Lambda)\to \Lambda$  by
$\Omega^u(x)=x^s$. The maps $\Omega^s$ and $\Omega^u$ are
$C^{\ell}$-smooth.

We now describe the scattering  map, following \cite{DelshamsLS08a}.  Assume that $W^u(\Lambda)$ has a transverse intersection with $W^s(\Lambda)$ along a $l$-dimensional homoclinic manifold $\Gamma$. The manifold $\Gamma$ consists of a $(l-1)$-dimensional family of trajectories asymptotic to $\Lambda$ in both forward and backwards time.
The transverse intersection of the hyperbolic invariant manifolds along $\Gamma$ means that  $\Gamma\subseteq W^u(\Lambda) \cap W^s(\Lambda)$
and, for each $x\in\Gamma$, we have
\begin{equation} \begin{split}\label{eq:dynamical channel}
T_xM=T_xW^u(\Lambda)+T_xW^s(\Lambda),\\
T_x\Gamma=T_xW^u(\Lambda)\cap T_xW^s(\Lambda).
\end{split} \end{equation}
Let us assume the additional condition that for each $x\in\Gamma$ we
have
\begin{equation} \label{eq:transverse foliation}
\begin{split}
T_xW^s(\Lambda)=T_xW^s(x^s)\oplus T_x(\Gamma),\\
T_xW^u(\Lambda)=T_xW^u(x^u)\oplus T_x(\Gamma),
\end{split}
\end{equation}
where $x^u,x^s$ are the uniquely defined points in $\Lambda$
corresponding to $x$.

The restrictions $\Omega^s_\Gamma,\Omega^u_\Gamma$ of
$\Omega^s,\Omega^u$, respectively, to $\Gamma$ are  local
$C^{\ell-1}$ - diffeomorphisms.  By restricting $\Gamma$ even further,  if necessary,  we can ensure that
$\Omega^s_\Gamma,\Omega^u_\Gamma$ are $C^{\ell-1}$-diffeomorphisms.
A homoclinic manifold $\Gamma$ for which the corresponding
restrictions of the wave maps are $C^{\ell-1}$-diffeomorphisms will
be referred as a homoclinic channel.

\begin{defn} \label{defn:scattering_map_flow}
Given a homoclinic channel $\Gamma$, the scattering map associated to
$\Gamma$ is the $C^{\ell-1}$-diffeomorphism
$S^\Gamma=\Omega_\Gamma^s\circ (\Omega_\Gamma^u)^{-1}$ defined on
the open subset $U^u:=\Omega_\Gamma^u(\Gamma)$ in $\Lambda$ to the
open subset $U^s:=\Omega_\Gamma^s(\Gamma)$ in $\Lambda$.
\end{defn}




\begin{prop}\label{prop:transversal}Assume that $T_1$ and $T_2$ are two invariant  submanifolds of complementary dimensions in $\Lambda$. Then  $W^u(T_1)$ has a transverse intersection with $W^s(T_2)$ in $M$ if and only if $S(T_1)$ has a transverse intersection with $T_2$ in $\Lambda$.\end{prop}

\section{Linearization of normally hyperbolic flows}
\label{section:linearization}
Let $M$ be a $C^r$-smooth, $m$-dimensional manifold (without boundary), with $r\geq 1$,
and  $\phi:M\times \mathbb{R}\to M$ a $C^r$-smooth flow
on $M$.
A submanifold  (possibly with boundary) $\Lambda$ of $M$ is said to
be a normally hyperbolic invariant manifold for $\phi^t$ if $\Lambda$ is invariant under $\phi^t$, there exists a splitting of the tangent bundle of $TM$ into sub-bundles
\[TM=E^u\oplus E^s\oplus T\Lambda,\]
that are invariant under $D\phi^t$ for all $t\in\mathbb{R}$, and there exist a constant $C>0$
and rates $0<\beta<\alpha$, such that for all $x\in\Lambda$ we have
\[\begin{split}\|D\phi^t(x)(v)\|\leq Ce^{-\alpha t}\|v\|  \textrm{ for all } t\geq 0, &\textrm{ if and only if } v\in E^s_x,\\
 \|D\phi^t(x)(v)\|\leq Ce^{\alpha t}\|v\|  \textrm{ for all } t\leq 0,  &\textrm{ if and only if }v\in E^u_x,\\
  \|D\phi^t(x)(v)\|\leq Ce^{\beta|t|
\|v\|} \textrm{ for all } t\in\mathbb{R}, &\textrm{ if and only if }v\in T_x\Lambda.
\end{split}\]

In \cite{HirschPS77} it is proved that in some neighborhood of $\Lambda$ the flow $\phi^t$ is conjugate with its linearization.
That is, there exists a neighborhood $\mathcal{U}$  of $\Lambda$ and a homeomorphism $h$ from $\mathcal{U}$ to some neighborhood $\mathcal{V}$
of the zero section  of the normal bundle to $\Lambda$ such that \[D\phi^t\circ h=h\circ \phi^t. \]

The homeomorphism $h$ defines a system of coordinates $(x_c,x_s,x_u)\in \Lambda\oplus E^s\oplus E^u$.  The flow written in these  coordinates takes the form
\[h(\phi^t(x))=(\phi^t(x_c), D\phi^t_{x_c}(x_s), D\phi^t_{x_c}(x_u)),\] for $x\in\mathcal{U}$ and $h(x)=(x_c,x_s,x_u)\in\mathcal{V}$.

\section{Correctly aligned windows.}\label{app:windows}

We follow \cite{GideaZ04a,GideaR03,GideaL06,GideaR12,DelshamsGR13}.

\begin{defn}
An $(m_1,m_2)$-window in an $m$-dimensional manifold $M$, where $m_1+m_2=m$, is a
compact subset $R$ of $M$ together with a $C^0$-parametrization given by a
homeomorphism $\rho$ from some open neighborhood $U_{R}$ of $[0,1]^{m_1}\times [0,1]^{m_2}\subseteq \mathbb{R}^{m_1}\times \mathbb{R}^{m_2}$ to an open subset of $M$, with $R=\rho([0,1]^{m_1}\times [0,1]^{m_2})$, and with a choice of an `exit set' \[R^{\rm
ex} =\rho \left(\partial[0,1]^{m_1}\times [0,1]^{m_2} \right )\]
and  of an `entry set'  \[R^{\rm en}
=\rho \left([0,1]^{m_1}\times
\partial[0,1]^{m_2}\right ).\]
\end{defn}
Let $f$ be a continuous
map on $M$ with $f(\textrm {im}(\rho_1))\subseteq \textrm
{im}(\rho_2)$.  Denote $f_\rho=\rho_2^{-1}\circ f\circ\rho_1$.
\begin{defn}\label{defn:corr}
Let  $R_1$ and $R_2$ be $(m_1,m_2)$-windows, and let $\rho_1$ and $\rho_2$ be the corresponding local parametrizations.
We say that $R_1$ is correctly aligned with
$R_2$ under $f$ if the following conditions are satisfied:
\begin{itemize}
\item[(i)] There exists a continuous homotopy $h:[0,1]\times
U_{R_1} \to {\mathbb R}^{m_1} \times {\mathbb R}^{m_2}$,
   such that the following conditions hold true
   \begin{eqnarray*}
      h_0&=&f_\rho, \\
      h([0,1],\partial[0,1]^{m_1}\times [0,1]^{m_2}) \cap  [0,1]^{m_1}\times [0,1]^{m_2} &=& \emptyset, \\
      h([0,1], [0,1]^{m_1}\times [0,1]^{m_2}) \cap [0,1]^{m_1}\times \partial[0,1]^{m_2} &=& \emptyset,
   \end{eqnarray*}
 and
\item[(ii)] There exists $y_0\in[0,1]^{m_2}$ such that
the map $A_{y_0}:\mathbb{R}^{m_1}\to\mathbb{R}^{m_1}$ defined by
$A_{y_0}(x)=\pi _{m_1}\left(h_{1}(x, y_0)\right )$ satisfies
\begin{eqnarray*}
A_{y_0}\left ( \partial[0,1]^{m_1}\right )\subseteq \mathbb
{R}^{m_1}\setminus [0,1]^{m_1},\\\deg({A_{y_0}},0)\neq 0,\end{eqnarray*}
where $\pi_{m_1}: \mathbb{R}^{m_1}\times \mathbb{R}^{m_2}\to
\mathbb{R}^{m_1}$ is the  projection onto the first
component, and $\deg$ is the Brouwer degree of the map $A_{y_0}$ at
$0$.
\end{itemize}
\end{defn}

\begin{thm}
\label{theorem:detorb} Let $\{R_i\}_{i\in\mathbb{Z}}$, be a collection of
$(m_1,m_2)$-windows in $M$,  and let $f_i$ be a collection of continuous maps on
$M$. If for each  $i\in\mathbb{Z}$, $R_i$ is correctly aligned with $R_{i+1}$ under $f_i$, then
there exists a point $p\in R_0$ such that
\[(f_{i}\circ \cdots\circ f_{0})(p)\in R_{i+1}, \textrm{ for all } i\in\mathbb{Z}. \]

Moreover, under the above conditions, and assuming that for some  $k>0$ we have  $R_{i}=R_{(i\,{\rm mod}\, k)}$  and $f_{i}=f_{(i\,{\rm mod}\, k)}$ for all
$i\in \mathbb{Z}$,
then there exists a
point $p$ as above that is periodic in the sense
\[(f_{k-1}\circ \cdots\circ f_{0})(p)=p.\]
\end{thm}

Assume that $f:M\to M$ is a diffeomorphism  on a manifold $M$, $\Lambda\subseteq M$ is an  $l$-dimensional normally hyperbolic invariant manifold, and $S:U\to V$ is a scattering map associated to some homoclinic channel $\Gamma$.
\begin{lem}\label{lem:shadowing2} Let $\{R_i,R'_i\}_{i\in\mathbb{Z}}$ be a bi-infinite
sequence of $l$-dimensional windows contained in  $\Lambda$. Assume that
the following properties hold for all $i\in\mathbb{Z}$:
\begin{itemize}
\item[(i)] $R_{i}\subseteq U$ and $R'_{i}\subseteq
V$.
\item[(ii)] $R_{i}$ is
correctly aligned with $R'_{i+1}$ under the scattering map
$S$.\item[(iii)] for  each pair $R'_{i+1},R_{i+1}$ and for each $L>0$ there exists $L'>L$  such that   $R'_{i+1}$ is correctly aligned
with   $R_{i+1}$ under the iterate $f_{\mid\Lambda}^{L'}$ of the restriction $f_{\mid\Lambda}$ of $f$ to $\Lambda$.
\end{itemize}
Fix any bi-infinite sequence of positive real numbers  $\{\eps_i\}_{i\in\mathbb{Z}}$.
Then there exist an orbit $(f^{n}(z))_{n\in\mathbb{Z}}$
of some point $z\in M$, an increasing sequence of integers
$(n_i)_{i\in\mathbb{Z}}$, and some sequences of positive integers $\{N_i\}_{i\in\mathbb{Z}}, \{K_i\}_{i\in\mathbb{Z}},
\{M_i\}_{i\in\mathbb{Z}}$,  such that, for all $i\in\mathbb{Z}$:
\[ \begin{split}
d(f^{n_i}(z),\Gamma)<\eps_i,\\
d(f^{n_i+N_{i+1}}(z), f_{\mid\Lambda}^{N_{i+1}}(R'_{i+1}))<\eps_{i+1},
\\d(f^{n_{i}-M_{i}}(z), f_{\mid\Lambda}^{-M_{i}}(R_{i}))<\eps_{i},\\
n_{i+1}=n_i+N_{i+1}+K_{i+1}+M_{i+1}.
\end{split}\]
\end{lem}

\section{Topological method for the diffusion problem.}\label{app:topmet}

In this section we recall the main result from \cite{GideaR12}.
Assume the following:
\begin{itemize}
\item [(C1)]  $M$ is a $n$-dimensional $C^r$-differentiable Riemannian manifold, and $f:M\to M$ is a $C^r$-smooth map, for some $r\geq 2$.
\item[(C2)] There exists a submanifold $\Lambda$ in $M$, diffeormorphic to an
annulus $\Lambda\simeq \mathbb{T}^1\times[0,1]$.  We assume that $f$ is  normally hyperbolic to $\Lambda$ in $M$.
Denote the dimensions of the stable and unstable manifolds of a point $x \in \Lambda$ by
$\text{dim}(W^s(x)) = n_s$ and $\text{dim}(W^u(x)) = n_u$.
Then, $n = 2 + n_s + n_u$.
\item[(C3)]  On $\Lambda$  there is a system of angle-action coordinates $(\phi,I)$, with
$\phi \in {\mathbb T}^1$ and $I \in [0,1]$.
The restriction $f|_{\Lambda}$ of $f$ to $\Lambda$ is a boundary component preserving,  area preserving, monotone twist map, with respect to the angle-action coordinates $(\phi,I)$.
\item [(C4)] The stable and unstable manifolds of $\Lambda$, $W^s(\Lambda)$ and $W^u(\Lambda)$, have a differentiably transverse intersection along a $2$-dimensional homoclinic channel $\Gamma$. We assume that the scattering map $S:U^-\to U^+$ associated to $\Gamma$  is well defined.
    \item [(C5)]
    There exists a bi-infinite sequence of Lipschitz primary invariant tori $\{T_i\}_{i\in\mathbb{Z}}$ in $\Lambda$, and
    a bi-infinite, increasing sequence of integers $\{i_k\}_{k\in\mathbb{Z}}$ with the following properties:
    \begin{itemize}
\item[(i)] Each torus $T_i$ intersects the domain $U^-$ and the range $U^+$ of the scattering map $S$ associated to $\Gamma$.
\item [(ii)] For each $i\in\{{i_{k}+1}, \ldots, {i_{k+1}-1}\}$, the image of $T_i\cap U^-$ under the scattering map $S$  is topologically transverse to  $T_{i+1}$.
\item [(iii)] For each torus $T_i$ with $i\in \{{i_{k}+2}, \ldots, {i_{k+1}-1}\}$, the  restriction of $f$ to $T_i$ is topologically transitive.
\item[(iv)] Each torus $T_i$ with $i\in \{{i_{k}+2}, \ldots, {i_{k+1}-1}\}$, can be
$C^0$-approximated from both sides by other primary invariant tori from
$\Lambda$.
\end{itemize}
We will refer to a finite sequence $\{T_i\}_{i=i_k+1,\ldots, i_{k+1}}$  as above as a transition chain of tori.
\item [(C6)]  The region in $\Lambda$ between $T_{i_k}$ and $T_{i_{k}+1}$ is a BZI.
\item[(C7)] Inside each region between $T_{i_k}$ and $T_{i_{k}+1}$   there is prescribed  a finite collection of Aubry-Mather sets $\{\Sigma _{\rho^k_1}, \Sigma _{\rho^k_2},\ldots, \Sigma _{\rho^k_{s_k}}\}$, where $s_k\geq 1$, and $\rho^k_s$ denotes the rotation number of $\Sigma_{\rho^k_s}$. These   Aubry-Mather sets  are assumed to be vertically ordered, relative to the $I$-coordinate on the annulus.
\end{itemize}

Instead of (C6)  we can consider the following condition:
\begin{itemize}
\item [(C6$'$)]  The region $\Lambda_k$ in $\Lambda$ between $T_{i_k}$ and $T_{i_{k}+1}$  contains finitely many invariant primary  tori $\{\Upsilon_{h^k_1}, \ldots , \Upsilon_{h^k_{l_k}}\}$, where $l_k\geq 1$, satisfying the following properties:
    \begin{itemize}
    \item [(i)] Each $\Upsilon_{h^k_j}$ falls in one of the following two cases:
    \begin{itemize}
    \item [(a)] $\Upsilon_{h^k_j}$ is an isolated invariant primary torus.
    \item [(b)] There exists a hyperbolic periodic orbit in $\Lambda$ such that its stable and unstable manifolds coincide.
    \end{itemize}
    \item [(ii)]  The invariant primary  tori $\{\Upsilon_{h^k_1}, \ldots , \Upsilon_{h^k_{l_k}}\}$ are vertically ordered, relative to the $I$-coordinate on the annulus.
    \item [(iii)] For each $\Upsilon_{h^k_j}$,  $j=1,\ldots, l_k$,
       the inverse image $S^{-1}(\Upsilon_{h^k_j}\cap U^+)$ forms with $\Upsilon_{h^k_j}$  a topological disk  $D_{h^k_j} \subseteq U^-$ below $\Upsilon_{h^k_j}$, such that $S(D_{h^k_j})\subseteq U^+$ is a topological disk  above $\Upsilon_{h^k_j}$, which is bounded by $\Upsilon_{h^k_j}$ and  $S(\Upsilon_{h^k_j}\cap U^-)$.
     \end{itemize}
\end{itemize}

\begin{thm}\label{thm:main1}
Let   $f:M\to M$  be a $C^r$-differentiable map, and let  $(T_i)_{i\in\mathbb{Z}}$ be a sequence of invariant primary
tori in $\Lambda$, satisfying the properties
(C1) -- (C6), or (C1)-(C5) and (C6\,$'$),  from above. Then for each  sequence
$(\epsilon_i)_{i\in\mathbb{Z}}$ of positive real numbers,   there exist a point $z\in M$
and a bi-infinite increasing sequence   of integers
$(N_i)_{i\in\mathbb{Z}}$  such that
\begin{eqnarray}\label{eqn:main11}  d(f^{N_i}(z), T_{i})<\epsilon _i, \textrm { for all }
i\in\mathbb{Z}.\end{eqnarray}

In addition, if  condition (C7) is assumed, and  some  positive integers $\{n^k_s\}_{s=1,\ldots, s_k}$, $k\in\mathbb{Z}$ are given, then there exist $z\in M$ and $(N_i)_{i\in\mathbb{Z}}$ as in \eqref{eqn:main11}, and positive integers $\{m^k_s\}_{s=1,\ldots, s_k}$, $k\in\mathbb{Z}$,  such that, for each $k$ and each $s\in\{1,\ldots, s_k\}$, we have
\begin{equation}\label{eqn:main22}\pi_\phi(f^j(w^k_s))<\pi_\phi(f^j(z))<\pi_\phi(f^j(\bar w^k_s)),\end{equation}
for some $w^k_s,\bar w^k_s\in\Sigma_{\omega^k_s}$ and for all $j$ with \[N_{i_k}+\sum_{t=0}^{s-1} n^k_t+\sum_{t=0}^{s-1}m^k_t\leq j\leq N_{i_k}+\sum_{t=0}^{s} n^k_t+\sum_{t=0}^{s-1}m^k_t.\]
\end{thm}

\bibliography{diffsym}
\bibliographystyle{plain}
\end{document}